\documentclass[11pt,letterpaper]{article}
\usepackage[letterpaper,margin=0.88in]{geometry}
\usepackage{amsmath,amssymb,mathtools,bm,booktabs,microtype}
\usepackage{pgfplots}
\usepackage{tikz}
\usepackage{cite}
\usepackage{hyperref}
\usepackage{graphicx}
\usepackage{float}
\usepackage[f]{esvect}

\allowdisplaybreaks

\hypersetup{hidelinks,pdfpagemode=UseNone,bookmarksopen=false}

\title{The $S$ and Double-Exponential Transformations with a Moment Expansion for Three-Centre Coulomb and Hybrid Integrals}
\author{Hassan Safouhi \\[1ex]
Campus Saint-Jean, University of Alberta,\\
Edmonton (AB), Canada\\
hsafouhi@ualberta.ca}
\date{}

\begin{document}

\maketitle

\vspace*{0.25cm}

\vspace*{0.25cm}
{\bf \large Abstract}. \hskip 0.10cm Three-centre two-electron Coulomb and hybrid integrals over exponential-type functions involve semi-infinite radial integrals containing oscillatory spherical Bessel functions. We investigate the combination of the $S$ transformation with Fourier-type double-exponential (DE) quadrature for their accurate and efficient numerical evaluation. The $S$ transformation transfers derivatives from the spherical-Bessel factor to the radial amplitude, converting the integral into a sine representation particularly well suited to Fourier-type DE quadrature. Numerical comparisons demonstrate that the combined $S+\mathrm{DE}$ approach can substantially reduce the number of quadrature nodes required to achieve high accuracy compared with direct DE. We further identify a geometry-dependent oscillation-to-damping ratio and establish a convergent moment expansion in the weakly oscillatory regime. Exploiting this complementary representation markedly enhances the numerical reliability of the algorithm in regimes where direct oscillatory quadrature becomes inefficient. The method is examined through semi-infinite radial integrals and complete H$_2$O and C$_2$H$_4$ molecular integrals over Slater-type functions, with comparisons against complete $S$-series evaluations and values reported in the literature.

\vspace*{0.25cm}

\section{Introduction}\label{introduction}
Three-centre two-electron Coulomb and hybrid integrals over exponential-type functions (ETFs) admit Fourier--$B$ representations in which numerical evaluation involves an outer Feynman-parameter integral and an inner semi-infinite integral containing an oscillatory spherical Bessel function~\cite{Safouhi11, Safouhi19}. The accurate evaluation of these radial integrals becomes particularly challenging for higher angular momenta and weakly damped oscillatory kernels. The $S$ transformation addresses this difficulty by transferring derivatives from the spherical-Bessel factor to the radial amplitude, converting the original integral into a Fourier-type sine integral~\cite{Safouhi19}. The equally spaced zeros of the sine function provide a computational advantage over the nonuniform finite zeros of higher-order spherical Bessel functions. Earlier work combined the $S$ transformation with Sidi's $\bar D$ transformation to accelerate the convergence of the transformed oscillatory integrals~\cite{Safouhi19}; the present work exploits a different numerical mechanism by combining the $S$ transformation with Fourier-type DE quadrature.

Earlier developments include the use of nonlinear transformations and properties of sine and Bessel functions for four-centre Coulomb integrals~\cite{Safouhi60, Safouhi38, Safouhi18, Safouhi19, Safouhi20}, high-accuracy numerical evaluations of three-centre nuclear-attraction and two-electron Coulomb and exchange integrals over Slater-type functions~\cite{Safouhi21}, and the treatment of two-centre Coulomb and hybrid integrals over Slater-type orbitals~\cite{Lesiuk-Moszynski-90-063318-14}. Recent developments confirm a continuing interest in the accurate and efficient evaluation of molecular integrals over ETFs. General representations for multicentre integrals over ETFs have recently been investigated~\cite{Yu-Si-Zhou-123-e27128-23}, while new numerical approaches have been developed for Slater-basis Coulomb and nuclear-attraction integrals in polyatomic calculations and for two-electron four-centre Coulomb integrals on modern GPU architectures~\cite{Stark-Meier-Hatch-Kammeraad-Dang-Zimmerman-47-e70291-26,Kaur-Chen-Luhana-Silberman-Ewert-Bunger-Segal-Greif-47-e70472-26}. These developments emphasize the importance of accurate and efficient methods for the evaluation of molecular integrals over ETFs.

The combination of the $S$ transformation with DE quadrature was previously applied to three-centre nuclear-attraction integrals~\cite{Safouhi55}. The two-electron Coulomb integrals considered here involve a richer radial structure. In particular, the balance between oscillation and molecular damping varies with the Feynman parameter, with the most pronounced effects occurring near the endpoints of the integration interval. We therefore examine the transformed radial integrals, their complex singularities, and the resulting behaviour of DE quadrature.

The combination of $S$ and DE exploits their complementary strengths. The $S$ transformation simplifies the oscillatory factor, while Fourier-type DE quadrature uses a mesh adapted to that structure and rapidly decreasing quadrature tails. In the strongly oscillatory test, direct DE requires 4695 quadrature nodes to reach an absolute error of $3.74\times10^{-15}$, whereas $S+\mathrm{DE}$ attains an absolute error of $4.36\times10^{-16}$ with only 107 nodes, a reduction by a factor greater than 43 in the number of quadrature nodes. The treatment of weakly oscillatory radial integrals by Berlu and Safouhi~\cite{Safouhi19} relied on neglecting sufficiently small higher-angular-momentum contributions and using Gauss--Laguerre quadrature for the limiting nonoscillatory case. The present work substantially extends this approach by introducing a convergent moment expansion that enables the systematic and accurate evaluation of these contributions, thereby improving the numerical reliability of the algorithm. Numerical tests on semi-infinite radial integrals and complete H$_2$O and C$_2$H$_4$ molecular integrals assess the complementary evaluation methods.

The remainder of the paper is organized as follows. Section~\ref{definitions-general-properties} presents the definitions and general properties. Section~\ref{molecular-reduction} gives the molecular reduction and the $S$ transformation. Section~\ref{double-exponential-quadrature} develops the DE formulations, while Section~\ref{convergence-analysis} examines their asymptotic behaviour and complex singularities. Section~\ref{geometry-moments} analyzes geometry-dependent damping and the moment expansion. Section~\ref{numerical-analysis} discusses the radial and molecular numerical results; Section~\ref{numerical-data} collects the corresponding tables and figures, followed by the conclusions.

\section{Definitions and General Properties}\label{definitions-general-properties}
The $B$ functions introduced by Steinborn and Filter are defined by~\cite{Steinborn-Filter-38-273-75, Filter-Steinborn-18-1-78}:
\begin{equation}
B_{n,l}^m(\zeta,\vec r)
=\frac{(\zeta r)^l}{2^{n+l}(n+l)!}
\widehat{k}_{n-\frac12}(\zeta r)
Y_l^m(\theta_{\vec r},\varphi_{\vec r}).
\label{b-function-definition}
\end{equation}
Here $\zeta>0$ is the orbital exponent, $(n, l, m)$ are the quantum numbers satisfying $n=1,2,\ldots$, $0\leq l\leq n-1$, and $-l\leq m\leq l$.

We use the surface spherical harmonics with the Condon--Shortley phase convention~\cite{Condon-Shortley-51}:
\begin{equation}
Y_l^m(\theta,\varphi)
=\mathrm i^{m+|m|}
\left(\frac{(2l+1)(l-|m|)!}{4\pi(l+|m|)!}\right)^{1/2}
P_l^{|m|}(\cos\theta)e^{\mathrm i m\varphi},
\label{spherical-harmonic-definition}
\end{equation}
where the associated Legendre polynomials.

The reduced Bessel functions are related to the modified Bessel functions of the second kind by~\cite{Watson-44}:
\begin{equation}
\widehat{k}_\nu(z)
=\sqrt{\frac{2}{\pi}}\,z^\nu K_\nu(z).
\label{khat-definition}
\end{equation}

They also satisfy the three-term recurrence relations~\cite{Abramowitz-Stegun-65}:
\begin{equation}
\begin{aligned}
\widehat{k}_{\nu+1}(z)
&=2\nu\widehat{k}_{\nu}(z)+z^2\widehat{k}_{\nu-1}(z).
\end{aligned}
\label{bessel-recurrences}
\end{equation}

The spherical Bessel functions satisfy the differential identity~\cite{Watson-44,Abramowitz-Stegun-65}:
\begin{equation}
 j_\lambda(z)=(-1)^\lambda z^\lambda
 \left(\frac{1}{z}\frac{d}{dz}\right)^\lambda\frac{\sin z}{z}.
\label{spherical-bessel-identity}
\end{equation}

They also satisfy the recurrence relation~\cite{Abramowitz-Stegun-65}:
\begin{equation}
x \, j_{l-1}(x) + x \, j_{l+1}(x) = (2l+1) \, j_{l}(x).
\label{EQRSHB}
\end{equation}

For nonnegative integer $\lambda$, the spherical-Bessel power series is~\cite{Watson-44,Abramowitz-Stegun-65,DLMF-26}:
\begin{equation}
j_\lambda(z)=\sum_{q=0}^{\infty}\frac{(-1)^q z^{\lambda+2q}}{2^q q!(2\lambda+2q+1)!!}.
\label{small-v-bessel-series}
\end{equation}

The Pochhammer symbol is given by~\cite{Abramowitz-Stegun-65}:
\begin{equation}
(\alpha)_n \,=\,
\left\{
  \begin{array}{l}
 (\alpha)_n \,= \, 1 \quad \textrm{if} \quad n = 0\\ 
 (\alpha)_n \, = \, \alpha \,(\alpha+1)\,(\alpha+2) \ldots (\alpha+n-1) \, = \,
\displaystyle \frac{\Gamma(\alpha+n)}{\Gamma(\alpha)}
\quad \textrm{if} \quad n \leq -\alpha\\ [0.2cm]
 (\alpha)_n \, = \, 0 \quad \textrm{if} \quad n \geq -\alpha+1,
  \end{array}
\right.
\label{pochhammer-definition}
\end{equation}
where $\Gamma$ stands for the Gamma function. 

The Gaunt coefficients are defined by:
\begin{equation}
\mathcal G(l_1m_1;l_2m_2;LM)
=\int_{S^2}\bigl[Y_{l_1}^{m_1}(\Omega)\bigr]^*Y_{l_2}^{m_2}(\Omega)Y_L^M(\Omega)\,d\Omega.
\label{gaunt}
\end{equation}

The normalized Slater-type functions (STFs) are defined by~\cite{Slater-42-33-32}:
\begin{equation}
\chi_{n,l}^{m}(\zeta,\vec r)
=\left(\frac{(2\zeta)^{2n+1}}{(2n)!}\right)^{1/2}r^{n-1}e^{-\zeta r}
Y_l^m(\theta_{\vec r},\varphi_{\vec r}),
\label{stf-definition}
\end{equation}
where $n=1,2,\ldots$, $0\leq l\leq n-1$, and $-l\leq m\leq l$.

STFs can be expressed as finite linear combinations of $B$ functions~\cite{Filter-Steinborn-18-1-78}:
\begin{equation}
\begin{aligned}
\chi_{n,l}^m(\zeta,\vec r)
={}&\frac{1}{\zeta^{n-1}}\left(\frac{(2\zeta)^{2n+1}}{(2n)!}\right)^{1/2}
\sum_{p=\widetilde p}^{n-l}
\frac{(-1)^{n-l-p}\,(n-l)!\,2^{2p+2l-n}(l+p)!}
{(2p-n+l)!\,(n-l-p)!}
B_{p,l}^m(\zeta,\vec r),
\end{aligned}
\label{stf-b-expansion}
\end{equation}
where the lower summation index is:
\begin{equation}
\widetilde p=\left\lceil\frac{n-l}{2}\right\rceil
=\begin{cases}
\dfrac{n-l}{2}, & n-l\ \text{even},\\
\dfrac{n-l+1}{2}, & n-l\ \text{odd}.
\end{cases}
\label{stf-b-lower-limit}
\end{equation}

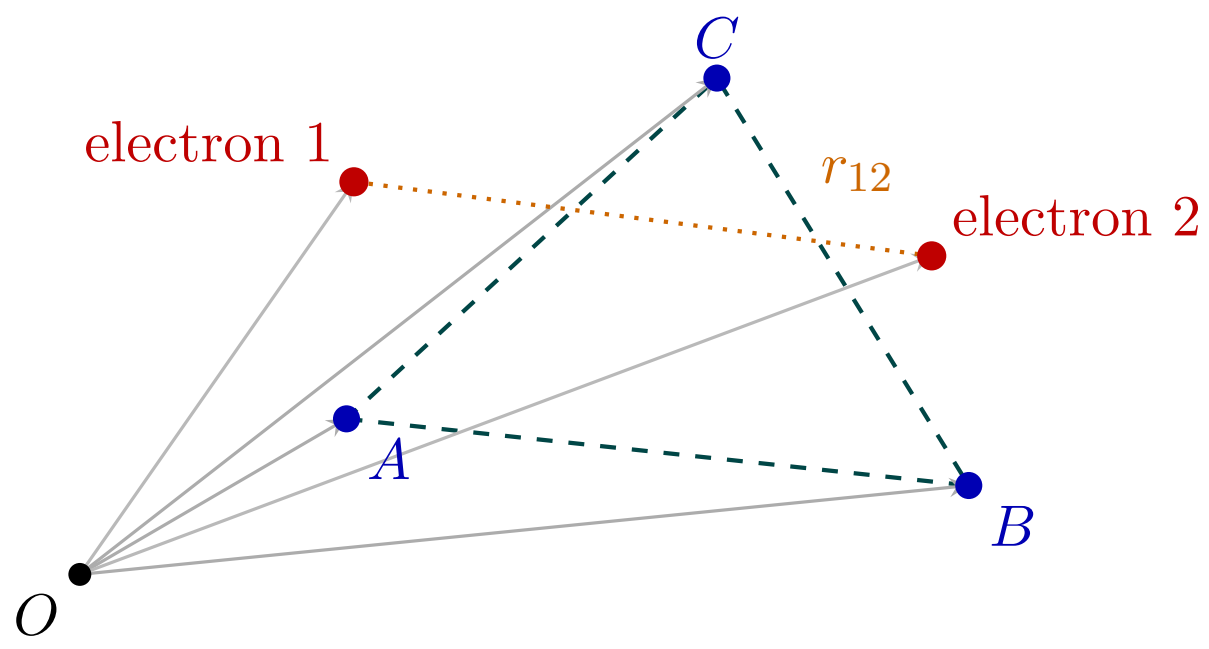
\begin{figure}[H]
\centering
\begin{tikzpicture}[x=1.0cm,y=1.0cm,>=stealth,scale=0.98]
\coordinate (O) at (0,0);
\coordinate (A) at (1.80,1.05);
\coordinate (B) at (6.0,0.60);
\coordinate (C) at (4.30,3.35);
\coordinate (eone) at (1.85,2.65);
\coordinate (etwo) at (5.75,2.15);
\draw[teal!55!black,dashed,thick] (A)--(B)--(C)--cycle;
\draw[->,gray!65,semithick] (O)--(A);
\draw[->,gray!65,semithick] (O)--(B);
\draw[->,gray!65,semithick] (O)--(C);
\draw[->,gray!55,semithick] (O)--(eone);
\draw[->,gray!55,semithick] (O)--(etwo);
\draw[dotted,orange!80!black,thick] (eone)--(etwo) node[pos=0.78,above=5pt,xshift=10pt] {$r_{12}$};
\fill[black] (O) circle(2.2pt) node[below left] {$O$};
\fill[blue!70!black] (A) circle(2.6pt) node[below right] {$A$};
\fill[blue!70!black] (B) circle(2.6pt) node[below right] {$B$};
\fill[blue!70!black] (C) circle(2.6pt) node[above] {$C$};
\fill[red!75!black] (eone) circle(2.8pt) node[above left] {electron 1};
\fill[red!75!black] (etwo) circle(2.8pt) node[above right] {electron 2};
\end{tikzpicture}
\caption{Two-dimensional projection of the three-centre geometry. The electron separation is $r_{12}=|\vec r_1-\vec r_2|$.}
\label{molecular-geometry}
\end{figure}

The three-centre two-electron Coulomb integral over $B$ functions is defined by:
\begin{equation}
\begin{aligned}
\mathcal K_{n_1l_1m_1,n_3l_3m_3}^{n_2l_2m_2,n_4l_4m_4}
={}&\iint_{\mathbb R^3\times\mathbb R^3}
\bigl[B_{n_1,l_1}^{m_1}(\zeta_1,\vec r_1-\vv{OA})\bigr]^*
B_{n_2,l_2}^{m_2}(\zeta_2,\vec r_1-\vv{OA})\\
&\quad\times\frac{1}{|\vec r_1-\vec r_2|}
\bigl[B_{n_3,l_3}^{m_3}(\zeta_3,\vec r_2-\vv{OB})\bigr]^*
B_{n_4,l_4}^{m_4}(\zeta_4,\vv r_2-\vec{OC})\,
 d\vec r_1\,d\vec r_2.
\end{aligned}
\label{molecular-integral-definition}
\end{equation}
The points $A$, $B$, and $C$ are the centres of the $B$ functions, and $\vec r_1$ and $\vec r_2$ denote the electron position vectors.
Figure~\ref{molecular-geometry} illustrates the three centres, the electron positions, and the fixed origin $O$.

By translation invariance, we move the coordinate origin from $O$ to $A$ and define:
\begin{equation}
\vec r_j^{\,\prime}=\vec r_j-\vv{OA}\quad (j=1,2),\quad
\vec R_3=\vv{OB}-\vv{OA},\quad
\vec R_4=\vv{OC}-\vv{OA},\quad
\vec R_{34}=\vec R_3-\vec R_4=\vv{OB}-\vv{OC}.
\label{translated-centres}
\end{equation}

Relabelling the translated electronic coordinates $\vec r_j^{\,\prime}$ as $\vec r_j$, the integral becomes:
\begin{equation}
\begin{aligned}
\mathcal K_{n_1l_1m_1,n_3l_3m_3}^{n_2l_2m_2,n_4l_4m_4}
={}&\iint_{\mathbb R^3\times\mathbb R^3}
\bigl[B_{n_1,l_1}^{m_1}(\zeta_1,\vec r_1)\bigr]^*
B_{n_2,l_2}^{m_2}(\zeta_2,\vec r_1)\\
&\quad\times\frac{1}{r_{12}}
\bigl[B_{n_3,l_3}^{m_3}(\zeta_3,\vec r_2-\vec R_3)\bigr]^*
B_{n_4,l_4}^{m_4}(\zeta_4,\vec r_2-\vec R_4)\,
 d\vec r_1\,d\vec r_2.
\end{aligned}
\label{molecular-integral-translated}
\end{equation}

The hybrid integral is obtained when $A$ and $B$ coincide, so that $\vec R_3=\vec 0$ in \eqref{molecular-integral-translated}. Three $B$ functions are then centred at $A$, while the fourth is centred at $C$:
\begin{equation}
\begin{aligned}
\mathcal H_{n_1l_1m_1,n_3l_3m_3}^{n_2l_2m_2,n_4l_4m_4}
={}&\iint_{\mathbb R^3\times\mathbb R^3}
\bigl[B_{n_1,l_1}^{m_1}(\zeta_1,\vec r_1)\bigr]^*
B_{n_2,l_2}^{m_2}(\zeta_2,\vec r_1)\frac{1}{r_{12}}
\bigl[B_{n_3,l_3}^{m_3}(\zeta_3,\vec r_2)\bigr]^*
B_{n_4,l_4}^{m_4}(\zeta_4,\vec r_2-\vec R_4)\,
 d\vec r_1\,d\vec r_2.
\end{aligned}
\label{hybrid-b-integral}
\end{equation}

Applying \eqref{stf-b-expansion} to each of the four orbital factors shows that the corresponding complete three-centre and hybrid integrals over STFs are finite linear combinations of the integrals over $B$ functions defined in \eqref{molecular-integral-translated} and \eqref{hybrid-b-integral}.

\section{Molecular Reduction and the $S$ Transformation}\label{molecular-reduction}
For the integral over $B$ functions defined in \eqref{molecular-integral-definition}, the established Fourier--$B$ reduction gives~\cite{Safouhi11, Safouhi19}:
\begingroup
\footnotesize
\begin{align}
\mathcal K_{n_1l_1m_1,n_3l_3m_3}^{n_2l_2m_2,n_4l_4m_4}={}&8(4\pi)^3\sqrt\pi\,
\frac{\zeta_1^{l_1}\zeta_2^{l_2}\zeta_3^{2n_3+l_3-1}\zeta_4^{2n_4+l_4-1}
(2l_3+1)!!(2l_4+1)!!(n_3+l_3+n_4+l_4+1)!}
{2^{l_1+l_2+1}\displaystyle\prod_{q=1}^{4}(n_q+l_q)!}\notag\\
&\times\sum_{l=|l_1-l_2|}^{l_1+l_2}\frac{(-i)^l\mathcal G(l_2m_2;l_1m_1;l,m_2-m_1)}{2^{2n_1+2n_2+l}}
\sum_{k=2}^{n_1+n_2}\sum_{i=k_1}^{k_2}
\frac{2^k(2n_1-i-1)!(2n_2-k+i-1)!\zeta_1^{i-1}\zeta_2^{k-i-1}}
{(i-1)!(n_1-i)!(k-i-1)!(n_2-k+i)!}\notag\\
&\times\sum_{l_4'=0}^{l_4}\sum_{m_4'=-l_4'}^{l_4'}\frac{i^{l_4+l_4'}(-1)^{l_4'}\mathcal G(l_4m_4;l_4-l_4',m_4-m_4';l_4'm_4')}{(2l_4'+1)!![2(l_4-l_4')+1]!!}\notag\\
&\times\sum_{l_3'=0}^{l_3}\sum_{m_3'=-l_3'}^{l_3'}\frac{i^{l_3+l_3'}\mathcal G(l_3m_3;l_3-l_3',m_3-m_3';l_3'm_3')}{(2l_3'+1)!![2(l_3-l_3')+1]!!}\notag\\
&\times\sum_{l'=0}^{l_3'+l_4'}\mathcal G(l_4'm_4';l_3'm_3';l',m_4'-m_3')R_{34}^{l'}Y_{l'}^{m_3'-m_4'}(\widehat{\vec R}_{34})\notag\\
&\times\sum_{l_{34}=0}^{(l_3-l_3')+(l_4-l_4')}\mathcal G(l_3-l_3',m_3-m_3';l_4-l_4',m_4-m_4';l_{34}m_{34})\notag\\
&\times\sum_{\lambda=0}^{l+l_{34}} i^\lambda\mathcal G(l,m_2-m_1;l_{34},m_{34};\lambda\mu)
\sum_{j=0}^{\Delta l}\binom{\Delta l}{j}\frac{(-1)^{(l_4'+l_3'+l')/2+j}}{2^{n_3+n_4+l_3+l_4-j+1}(n_3+n_4+l_3+l_4-j+1)!}\notag\\
&\times\frac{\zeta_s^{n_k-l-1}\Gamma(k+l_1+l_2+l+1)}{\Gamma(l+3/2)}
\sum_{r=0}^{\eta'}(-1)^r\frac{(\eta/2)_r((\eta+1)/2)_r}{(l+3/2)_r r!\zeta_s^{2r}}\notag\\
&\times\int_0^1s^{n_3+l_3+l_4-l_4'}(1-s)^{n_4+l_4+l_3-l_3'}Y_\lambda^\mu(\widehat{\vec v}(s))\left[\int_0^\infty \frac{x^{n_x}}{(\zeta_s^2+x^2)^{n_k}}
\frac{\widehat{k}_\nu(R_{34}\gamma(s,x))}{\gamma(s,x)^{n_\gamma}}j_\lambda(vx)\,dx \right]\,ds.
\label{auto-034}
\end{align}
\endgroup

The parameters $\zeta_1$, $\zeta_2$, $\zeta_3$, and $\zeta_4$ are the positive orbital exponents of the four respective $B$ functions. The quantities involved in the complete molecular expression are:
\begin{align}
\zeta_s&=\zeta_1+\zeta_2,
\quad k_1=\max(1,k-n_2)
\quad\text{and}\quad k_2=\min(n_1,k-1),
\notag\\[0.35em]
n_x&=l_3-l_3'+l_4-l_4'+2r+l
\quad\text{and}\quad n_k=k+l_1+l_2,
\notag\\[0.35em]
n_\gamma&=2(n_3+l_3+n_4+l_4)-(l_3'+l_4')-l'+1,
\notag\\[0.35em]
\Delta l&=\frac{l_3+l_4-l'}{2}
\quad\text{with}\quad \eta=l-k-l_1-l_2+1,
\notag\\[0.35em]
\eta'&=
\begin{cases}
-\eta/2 & \text{if} \quad \eta \text{ is even},\\
-(\eta+1)/2 & \text{if} \quad \eta \text{ is odd},
\end{cases}
\notag\\[0.35em]
m_{34}&=(m_3-m_3')-(m_4-m_4')
\quad\text{with}\quad \mu=(m_2-m_1)-m_{34},
\notag\\[0.35em]
\nu&=n_3+n_4+l_3+l_4-l'-j+\tfrac12,
\notag\\[0.35em]
\gamma(s,x)^2&=(1-s)\zeta_4^2+s\zeta_3^2+s(1-s)x^2,\notag\\[0.35em]
\vec v(s)&=(1-s)\vec R_{34}+\vec R_4
\quad\text{with}\quad v(s)=|\vec v(s)|.
\label{parameter-set}
\end{align}

In the following, we write $v$ for the $s$-dependent magnitude $v(s)$.

For hybrid integrals $\mathcal H_{n_1l_1m_1,n_3l_3m_3}^{n_2l_2m_2,n_4l_4m_4}$, the analytical representation follows directly from the three-centre expression by setting $\vec R_3=\vec 0$. Consequently, $\vec v(s)=s\vec R_4$ and $v(s)=sR_4$.

In the representation~\eqref{auto-034}, the semi-infinite radial integral is:
\begin{equation}
\mathcal K(s)=\int_0^\infty \frac{x^{n_x}}{(\zeta_s^2+x^2)^{n_k}}
\frac{\widehat{k}_\nu(R_{34}\gamma(s,x))}{\gamma(s,x)^{n_\gamma}}j_\lambda(vx)\,dx.
\label{auto-008}
\end{equation}

The corresponding semi-infinite radial integral for the hybrid case is denoted by $\mathcal H(s)$, and is obtained from $\mathcal K(s)$ by by setting $\vec R_3=\vec 0$ and $\vec v(s)=s\vec R_4$ .

For $v>0$, substituting $z=vx$ into \eqref{spherical-bessel-identity} gives~\cite{Safouhi15,Safouhi19}:
\begin{equation}
\mathcal K(s)=\frac{(-1)^\lambda}{v^{\lambda+1}}\int_0^\infty \frac{x^{n_x+\lambda}}{(\zeta_s^2+x^2)^{n_k}}\frac{\widehat{k}_\nu(R_{34}\gamma(s,x))}{\gamma(s,x)^{n_\gamma}}
\left(\frac{1}{x}\frac{d}{dx}\right)^\lambda\!\left[\frac{\sin(vx)}{x}\right]\,dx.
\label{s-before-parts}
\end{equation}

The $S$ transformation relies on a reformulated integration-by-parts procedure with respect to the measure $x\,dx$~\cite{Safouhi15,Safouhi19}. Repeated application transfers the derivatives from the spherical Bessel factor to the radial factor. The boundary terms vanish under the following endpoint conditions.

As $x\to0^+$, the differentiated radial factor satisfies:
\begin{equation}
\left(\frac{1}{x}\frac{d}{dx}\right)^\ell\!\left[\frac{x^{n_x+\lambda-1}}{(\zeta_s^2+x^2)^{n_k}}\frac{\widehat{k}_\nu(R_{34}\gamma(s,x))}{\gamma(s,x)^{n_\gamma}}\right]=O\!\left(x^{n_x+\lambda-1-2\ell}\right),\qquad 0\leq\ell\leq\lambda-1.
\label{S-origin}
\end{equation}

The Gaunt selection rules in \eqref{auto-034} give:
\begin{equation}
l_{34}\leq(l_3-l_3')+(l_4-l_4')  \quad \text{and} \quad \lambda\leq l+l_{34}.
\end{equation}

Since $r$ is the nonnegative summation index in the finite hypergeometric expansion, the angular-momentum coupling conditions imply:
\begin{equation}
\lambda\leq l+l_{34}\leq l+(l_3-l_3')+(l_4-l_4')\leq n_x=l+(l_3-l_3')+(l_4-l_4')+2r.
\label{nx-lambda}
\end{equation}

The smallest exponent in \eqref{S-origin} is:
\begin{equation}
n_x-\lambda+1\geq1.
\end{equation}
The boundary terms at the origin arising from the reformulated integration by parts therefore vanish.

For fixed $0<s<1$ and $R_{34}>0$, the Macdonald factor gives:
\begin{equation}
\left(\frac{1}{x}\frac{d}{dx}\right)^\ell\!\left[\frac{x^{n_x+\lambda-1}}{(\zeta_s^2+x^2)^{n_k}}\frac{\widehat{k}_\nu(R_{34}\gamma(s,x))}{\gamma(s,x)^{n_\gamma}}\right]=O\!\left(x^{\eta_\ell}e^{-R_{34}\sqrt{s(1-s)}x}\right) \quad \text{as} \quad x\to\infty,
\label{S-infinity}
\end{equation}
for finite algebraic exponents $\eta_\ell$. The corresponding boundary terms at infinity also vanish.

After $\lambda$ applications of the reformulated integration-by-parts identity, the factor $(-1)^\lambda$ in \eqref{s-before-parts} cancels. Thus:
\begin{equation}
\mathcal K(s)=\frac{1}{v^{\lambda+1}}\int_0^\infty \left(\frac{1}{x}\frac{d}{dx}\right)^\lambda\!\left[\frac{x^{n_x+\lambda-1}}{(\zeta_s^2+x^2)^{n_k}}\frac{\widehat{k}_\nu(R_{34}\gamma(s,x))}{\gamma(s,x)^{n_\gamma}}\right]\sin(vx)\,dx.
\label{auto-009}
\end{equation}

For comparison with numerical evaluations, let $G_{\lambda,s}(x)$ denote the differentiated radial factor in \eqref{auto-009}. Integrating over successive half-periods of the sine function gives the $S$-series representation:
\begin{equation}
\mathcal K(s)=\frac{1}{v^{\lambda+1}}\sum_{n=0}^{\infty}\int_{n\pi/v}^{(n+1)\pi/v}G_{\lambda,s}(x)\sin(vx)\,dx.
\label{s-series}
\end{equation}

The complete $S$-series values are obtained by performing the outer Feynman-parameter integration and the finite summation over $B$-function integrals with the coefficients of the STF expansions.

\section{Double-exponential quadrature and the two numerical formulations}\label{double-exponential-quadrature}

The two quadrature formulations operate on different radial representations:
\begin{equation}
\mathcal O(vx)=
\begin{cases}
j_\lambda(vx) & \text{for} \quad \text{direct DE},\\
\sin(vx) &  \text{for} \quad S+\mathrm{DE}.
\end{cases}
\label{two-kernels}
\end{equation}
The following DE formulations evaluate the spherical-Bessel and sine representations. A complementary moment evaluation for small oscillation-to-damping ratios is developed in Section~\ref{geometry-moments}.

Related numerical developments include error analyses of transformed trapezoidal rules for Fourier-type integrals~\cite{Denich-Novati-198-160-24}, Sinc quadrature for Hankel transforms~\cite{Denich-Novati-100-23-24}, and Clenshaw--Curtis product-integration methods for highly oscillatory or exponentially decaying integrals~\cite{Dominguez-26}.

The DE transformation and trapezoidal discretization follow \cite{Takahasi-Mori-9-721-74, Mori-12-119-85, Ooura-Mori-38-353-91, Ooura-Mori-112-229-99, Mori-Sugihara-127-287-01, Mori-41-897-05}. We consider:
\begin{equation}
 I=\int_0^\infty F(x)\mathcal O(vx)\,dx,
\label{auto-018}
\end{equation}
and set $x=M\phi(t)$:
\begin{equation}
 I=M\int_{-\infty}^{\infty}F(M\phi(t))\mathcal O(vM\phi(t))\phi'(t)\,dt.
\label{auto-019}
\end{equation}

The infinite and truncated trapezoidal approximations are:
\begin{equation}
 I_h=Mh\sum_{k=-\infty}^{\infty}F(M\phi(kh))\mathcal O(vM\phi(kh))\phi'(kh),
\label{auto-020}
\end{equation}
and:
\begin{equation}
 I_{h,N_-,N_+}=Mh\sum_{k=-N_-}^{N_+}F(M\phi(kh))\mathcal O(vM\phi(kh))\phi'(kh).
\label{auto-021}
\end{equation}

The separate indices $N_-$ and $N_+$ allow asymmetric transformed tails. The error decomposition is:
\begin{equation}
 I-I_{h,N_-,N_+}=(I-I_h)+(I_h-I_{h,N_-,N_+}),
\label{auto-022}
\end{equation}
The first term is the discretization error, and the second is the truncation error. Section~\ref{convergence-analysis} examines the analyticity and asymptotic properties relevant to the discretization and decay of the sampled terms.

For the transformed sine integral, the zeros of the oscillatory factor occur exactly at:
\begin{equation}
 x_n=\frac{n\pi}{v} \quad \text{for} \quad n=0,1,2,\ldots .
\label{sine-zeros}
\end{equation}

The Fourier-type DE construction exploits these regularly spaced zeros by choosing:
\begin{equation}
 Mh=\frac{\pi}{v}.
\label{auto-023}
\end{equation}

Two Fourier-type mappings are considered for the $S$-transformed sine integral. The first is the Ooura--Mori mapping~\cite{Ooura-Mori-38-353-91, Mori-41-897-05}:
\begin{equation}
\phi_1(t)=\frac{t}{1-\exp[-K\sinh t]} \quad \text{with} \quad K=6.
\label{auto-024}
\end{equation}
Following Ooura and Mori~\cite{Ooura-Mori-38-353-91,Mori-41-897-05}, we choose $K=6$. This value is close to $2\pi$ while remaining below the threshold beyond which the complex singularities of the mapping approach the real axis. It therefore provides rapid double-exponential decay without reducing the analyticity strip of the mapping. 
The numerical results in Table~\ref{K-sensitivity} identify $K=6$ as the optimal choice among the values tested, consistently achieving near-machine-precision accuracy across all four representative kernels while maintaining a low number of quadrature nodes.

The second mapping is the improved Fourier-type DE mapping of Ooura and Mori~\cite{Ooura-Mori-112-229-99}:
\begin{equation}
\phi_2(t)=\frac{t}{1-\exp[-Q_2(t)]}
\quad \text{with} \quad Q_2(t)=2t+\alpha(1-e^{-t})+\beta(e^t-1),
\label{phi2}
\end{equation}
with:
\begin{equation}
\beta=\frac14
\quad\text{and}\quad
\alpha=\frac{\beta}{\sqrt{1+M\log(1+M)/(4\pi)}}.
\label{phi2-parameters}
\end{equation}
Both mappings satisfy $\phi_j(t)\to0$ double exponentially as $t\to-\infty$ and $\phi_j(t)-t\to0$ double exponentially as $t\to+\infty$. Consequently, \eqref{auto-023} gives for $j=1,2$:
\begin{equation}
vM\phi_j(kh)=k\pi+vM[\phi_j(kh)-kh]\longrightarrow k\pi,
\label{zero-alignment-limit}
\end{equation}
so the positive DE nodes approach the zeros of the sine function for both mappings.

For $j_\lambda(vx)$ the zeros are not equally spaced at finite index. Since its asymptotic period is $2\pi/v$, direct DE is supplemented by the conservative phase-resolution criterion:
\begin{equation}
vMh\leq\frac{\pi}{2},
\label{direct-phase-resolution}
\end{equation}
which gives asymptotically at least four nodes per oscillation and is used to ensure adequate resolution of the oscillations. Equations \eqref{auto-023} and \eqref{direct-phase-resolution} therefore reflect the different numerical mechanisms of asymptotic alignment with the sine zeros and direct phase resolution.

\section{Asymptotic Behaviour and Complex Singularities}\label{convergence-analysis}

The two Fourier-type DE mappings satisfy $\phi_j(t)\to0$ as $t\to-\infty$ and $\phi_j(t)-t\to0$ as $t\to+\infty$. The second property makes the positive quadrature nodes approach the zeros of the sine factor in the $S$-transformed integral.

For $\phi_1$, the relation $\sinh t\sim e^t/2$ as $t\to+\infty$ gives:
\begin{equation}
\phi_1(t)-t=\frac{t e^{-K\sinh t}}{1-e^{-K\sinh t}}\sim t\exp\!\left(-\frac{K}{2}e^t\right),\qquad t\to+\infty.
\label{auto-026}
\end{equation}

At the opposite limit, $\sinh t\sim-e^{-t}/2$, and hence:
\begin{equation}
\phi_1(t)\sim(-t)\exp\!\left(-\frac{K}{2}e^{-t}\right),\qquad t\to-\infty.
\label{phi1-negative-scale}
\end{equation}

For $\phi_2$, the exponent $Q_2(t)$ in \eqref{phi2} yields:
\begin{equation}
\phi_2(t)-t\sim t\exp\!\left(-\beta e^t-2t-\alpha+\beta\right),\qquad t\to+\infty,
\label{phi2-positive}
\end{equation}
and:
\begin{equation}
\phi_2(t)\sim(-t)\exp\!\left(-\alpha e^{-t}+2t+\alpha-\beta\right),\qquad t\to-\infty.
\label{phi2-negative}
\end{equation}

Thus $\phi_j(t)$ approaches zero double exponentially at the negative end of the real axis, whereas $\phi_j(t)-t$ decreases double exponentially at the positive end.

Using $vMh=\pi$, define the deviation of the sampled phase from a sine zero by:
\begin{equation}
vM\phi_j(kh)=k\pi+\delta_{j,k},\qquad \delta_{j,k}=vM[\phi_j(kh)-kh].
\label{zero-deviation}
\end{equation}

The positive-side asymptotics then give:
\begin{equation}
\begin{aligned}
\delta_{1,k}&\sim k\pi\exp\!\left(-\frac{K}{2}e^{kh}\right),\\
\delta_{2,k}&\sim k\pi\exp\!\left(-\beta e^{kh}-2kh-\alpha+\beta\right),
\end{aligned}\qquad k\to+\infty.
\label{auto-027}
\end{equation}

Consequently:
\begin{equation}
\sin[vM\phi_j(kh)]=(-1)^k\sin\delta_{j,k}\sim(-1)^k\delta_{j,k},\qquad k\to+\infty.
\label{auto-028}
\end{equation}

The sampled sine factor therefore decreases double exponentially for both maps. This accounts for the rapid decrease of the positive summation terms in the $S+\mathrm{DE}$ quadrature.

We next identify the finite complex singularities of the radial integrand. For fixed $0<s<1$, define:
\begin{equation}
a_s^2=\frac{(1-s)\zeta_4^2+s\zeta_3^2}{s(1-s)},\qquad \gamma(s,x)^2=s(1-s)(x^2+a_s^2).
\label{gamma-singularities}
\end{equation}

The radial factor in \eqref{auto-008} is locally $x^{n_x}A_s(x^2)$, where $A_s$ is analytic near zero because $\zeta_s>0$ and $\gamma(s,0)>0$. The differentiated factor in \eqref{auto-009} therefore satisfies:
\begin{equation}
\left(\frac{1}{x}\frac{d}{dx}\right)^\lambda\!\left[x^{n_x+\lambda-1}A_s(x^2)\right]=x^{n_x-\lambda-1}B_s(x^2),
\label{origin-derivative}
\end{equation}
where $B_s$ is analytic near zero. Multiplying by the sine factor gives:
\begin{equation}
\sin(vx) \left(\frac{1}{x}\frac{d}{dx}\right)^\lambda\!\Bigg[x^{n_x+\lambda-1}A_s(x^2)\Bigg]=x^{n_x-\lambda}B_s(x^2)\frac{\sin(vx)}{x}.
\label{origin-regularity}
\end{equation}

Since $n_x\geq\lambda$ by \eqref{nx-lambda} and $\sin(vx)/x$ is analytic at $x=0$, the complete $S$-transformed integrand is analytic at the origin.

Away from the origin, the radial factors $(\zeta_s^2+x^2)^{-n_k}$ and $\gamma(s,x)$ give the possible finite singularity locations:
\begin{equation}
\Sigma_x(s)=\{\pm i\zeta_s,\ \pm ia_s\}.
\label{molecular-singularity-set}
\end{equation}

The half-integer reduced Bessel factor and the operator $D^\lambda$ may change the local nature or order of these singularities, but introduce no additional nonzero finite singularity locations.

The first DE map has singularities where its denominator vanishes:
\begin{equation}
\Sigma_{\phi_1}=\{t\in\mathbb C:K\sinh t=2\pi i m,\ m\in\mathbb Z\}\setminus\{0\}.
\label{phi1-singularities}
\end{equation}

Since the point $t=0$ is removable, the distance from the real axis to the nearest nonremovable singularities is:
\begin{equation}
\operatorname{dist}(\mathbb R,\Sigma_{\phi_1})=
\begin{cases}
\dfrac{\pi}{2},&0<K\leq2\pi,\\[1mm]
\arcsin\!\left(\dfrac{2\pi}{K}\right),&K>2\pi.
\end{cases}
\label{phi1-strip-distance}
\end{equation}

For $\phi_2$, the corresponding singularities satisfy:
\begin{equation}
\Sigma_{\phi_2}=\{t\in\mathbb C:Q_2(t)=2\pi i m,\ m\in\mathbb Z\}\setminus\{0\}.
\label{phi2-singularities}
\end{equation}

Their locations depend on $M$ through $\alpha$ in \eqref{phi2-parameters}.

The inverse images of the radial singularities under $x=M\phi_j(t)$ impose an additional restriction. Define:
\begin{equation}
d_{\mathrm{mol},j}(s,M)=\operatorname{dist}\!\left(\mathbb R,\{t:M\phi_j(t)\in\Sigma_x(s)\}\right),\qquad j=1,2.
\label{molecular-pullback-distance}
\end{equation}

An admissible strip half-width $d_j$ must satisfy:
\begin{equation}
0<d_j<\min\!\left\{\frac{\pi}{2},\operatorname{dist}(\mathbb R,\Sigma_{\phi_j}),d_{\mathrm{mol},j}(s,M)\right\},\qquad j=1,2.
\label{strip-bound}
\end{equation}

The restriction $d_j<\pi/2$ ensures that the asymptotic estimates remain valid as the real part of $t$ tends to $\pm\infty$ within the analyticity strip. For fixed $0<s<1$ and $R_{34}>0$, if the transformed integrand is analytic in this strip and satisfies the required horizontal-line integrability conditions, the infinite trapezoidal approximation has discretization error~\cite{Mori-Sugihara-127-287-01,Sugihara-72-767-02}:
\begin{equation}
|I-I_h|=O\!\left(e^{-2\pi d_j/h}\right),\qquad h\to0^+.
\label{strip-error}
\end{equation}

\section{Geometry-Dependent Radial Behaviour and Moment Expansion}\label{geometry-moments}

The endpoint and large-$x$ structure follows directly from \eqref{parameter-set}. For fixed $x$:
\begin{equation}
\gamma(s,x)=\zeta_4+\frac{s}{2\zeta_4}(\zeta_3^2-\zeta_4^2+x^2)+O(s^2) \quad \text{as} \quad s\to0^+.
\label{auto-012}
\end{equation}

At the opposite endpoint, we set $\varepsilon=1-s$, so that $s\to1^-$ corresponds to $\varepsilon\to0^+$. For fixed $x$, we obtain:
\begin{equation}
\gamma(1-\varepsilon,x)=\zeta_3+\frac{\varepsilon}{2\zeta_3}(\zeta_4^2-\zeta_3^2+x^2)+O(\varepsilon^2)  \quad \text{as} \quad \varepsilon\to0^+.
\label{auto-013}
\end{equation}

For fixed $0<s<1$ and $x\to\infty$:
\begin{equation}
\gamma(s,x)=\sqrt{s(1-s)}\,x\left[1+\frac{(1-s)\zeta_4^2+s\zeta_3^2}{2s(1-s)x^2}+O(x^{-4})\right].
\label{auto-014}
\end{equation}

Using the standard large-argument Macdonald-function expansion \cite{DLMF-26}:
\begin{equation}
K_\nu(z)\sim\sqrt{\frac{\pi}{2z}}e^{-z}\left[1+\frac{4\nu^2-1}{8z}+O(z^{-2})\right],
\label{auto-015}
\end{equation}
this gives the dominant molecular damping scale:
\begin{equation}
\exp[-R_{34}\sqrt{s(1-s)}\,x].
\label{auto-016}
\end{equation}

Its decay rate collapses at the Feynman endpoints.

For fixed $0<s<1$ and $R_{34}>0$, the decay rate is therefore:
\begin{equation}
\alpha_s=R_{34}\sqrt{s(1-s)}.
\label{alphas}
\end{equation}

The corresponding dimensionless ratio is:
\begin{equation}
\rho(s)=\frac{v}{R_{34}\sqrt{s(1-s)}}.
\label{auto-017}
\end{equation}

For hybrid integrals, $v=sR_4$, which gives:
\begin{equation}
\rho(s)=\sqrt{\frac{s}{1-s}} \quad \Longrightarrow  \quad \rho(s)<1   \quad \text{for} \quad  0<s<1/2.
\end{equation}

In the case of three-centre integrals, $v$ does not generally vanish at either endpoint. It can become small in the interior when $\vec R_3$ and $\vec R_4$ partly cancel, and vanishes for an antiparallel collinear assignment at $s=R_3/(R_3+R_4)$. Hence the numerical regimes depend on the assignment of the centres as well as on their geometry.

The ratio $\rho(s)$ compares oscillation with exponential damping. The condition $\rho(s)<1$ also supplies a sufficient convergence domain for the nonoscillatory moment expansion.

Using the spherical-Bessel power series \eqref{small-v-bessel-series}, its absolute sum satisfies, for fixed $\lambda$:
\begin{equation}
\sum_{q=0}^{\infty}
\frac{|v|^{\lambda+2q}x^{\lambda+2q}}
{2^q q!(2\lambda+2q+1)!!}
\leq C_\lambda (|v|x)^\lambda e^{|v|x},
\label{small-v-bessel-majorant}
\end{equation}
with a constant $C_\lambda>0$ depending only on $\lambda$.

Together with the Macdonald decay, this gives the absolute majorant:
\begin{equation}
\sum_{q=0}^{\infty}\left|
\frac{v^{\lambda+2q}x^{n_x+\lambda+2q}}
{2^q q!(2\lambda+2q+1)!!}
\frac{\widehat{k}_\nu(R_{34}\gamma(s,x))}
{(\zeta_s^2+x^2)^{n_k}\gamma(s,x)^{n_\gamma}}
\right|
=O\!\left(x^{\chi}e^{-(\alpha_s-|v|)x}\right),
\label{small-v-majorant}
\end{equation}
for a finite algebraic exponent $\chi$.
Thus $\rho(s)<1$ is a sufficient condition for termwise integration in the radial integral~\eqref{auto-008}.

The bound in \eqref{small-v-majorant} permits the limit of the spherical-Bessel partial sums to be passed through the integral, giving the moment representation:
\begin{equation}
\mathcal K(s)=\sum_{q=0}^{\infty}\frac{(-1)^q v^{\lambda+2q}}{2^q q!(2\lambda+2q+1)!!}\,\mathcal M_q(s)
\quad \text{with} \quad
\mathcal M_q(s)=\int_0^\infty \frac{x^{n_x+\lambda+2q}}{(\zeta_s^2+x^2)^{n_k}}\frac{\widehat{k}_\nu(R_{34}\gamma(s,x))}{\gamma(s,x)^{n_\gamma}}\,dx.
\label{small-v-moments}
\end{equation}

At fixed $0<s<1$, the leading term makes the regular small-$v$ limit explicit:
\begin{equation}
\mathcal K(s)=\frac{v^\lambda}{(2\lambda+1)!!}\mathcal M_0(s)+O(v^{\lambda+2})
\quad \text{as} \quad v\to0^+.
\label{small-v-leading}
\end{equation}
The limiting regime $v\to0$ was previously considered by Berlu and Safouhi~\cite{Safouhi19}. In that treatment, contributions associated with $\lambda>0$ were regarded as negligible for sufficiently small $v$, while the $\lambda=0$ limit was reduced to a nonoscillatory semi-infinite integral evaluated by Gauss--Laguerre quadrature. The present development extends this treatment by establishing the leading $v^\lambda$ dependence and a convergent moment representation under the sufficient condition $\rho(s)<1$. The higher-angular-momentum contributions can therefore be evaluated systematically for both hybrid and three-centre integrals, rather than neglected solely because $v$ is small.

Near a geometric zero of $v$, the moment representation avoids the prefactor $v^{-(\lambda+1)}$ in the $S$-transformed expression. At $v=0$, the original integral is interpreted through the regular limit.

This fixed-$s$ argument does not furnish a uniform endpoint estimate because $\alpha_s\to0$ as $s\to0^+$ or $1^-$. Nevertheless, the hybrid geometry satisfies $\rho(s)\to0$ as $s\to0^+$, making the moment representation particularly effective near the lower endpoint. The significance of this development extends beyond the derivation of an alternative series representation. It replaces the earlier small-$v$ approximation, in which higher-angular-momentum contributions were neglected, with a convergent and systematically computable expansion that retains these contributions and avoids the inverse powers of $v$ appearing in the $S$-transformed formulation. Implemented as a separate small-$\rho$ evaluation branch for both hybrid and three-centre integrals, this approach substantially improves the numerical robustness of the algorithm in weakly oscillatory regimes. As demonstrated in~\ref{small-v-diagnostics}, the moment expansion achieves errors many orders of magnitude smaller than those of an unmodified fixed-grid $S+\mathrm{DE}$ calculation in the tested small-$\rho$ configurations.

\section{Numerical Analysis and Discussion}\label{numerical-analysis}

The numerical results in Section~\ref{numerical-data} were generated in double-precision. The $\phi_1$ calculations use $K=6$, while $\phi_2$ uses \eqref{phi2-parameters}. Both quadrature sums are truncated when their terms become smaller than $10^{-15}$ in absolute value. The outer $s$-integration uses a 48-point Gauss--Legendre quadrature. Except for Table~\ref{small-v-diagnostics}, which uses direct integration of the original spherical-Bessel radial integrands as its reference, the absolute errors reported in the remaining radial tables are computed relative to values obtained from the complete $S$-series~\eqref{s-series}. The complete molecular-integral tables likewise use the corresponding $S$-series evaluations as reference values. 

The first numerical results concern the semi-infinite integrals $\mathcal H(s)$ and $\mathcal K(s)$ that occur in the Fourier--$B$ representation of the H$_2$O molecular integrals. We use $\langle 1s^{\mathrm{O}}1s^{\mathrm{O}}\mid1s^{\mathrm{O}}1s^{\mathrm{H}_1}\rangle$ for the hybrid case and $\langle 1s^{\mathrm{O}}1s^{\mathrm{O}}\mid1s^{\mathrm{H}_1}1s^{\mathrm{H}_2}\rangle$ for the three-centre case.

\begin{figure}[H]
\centering
\includegraphics[width=0.6\linewidth]{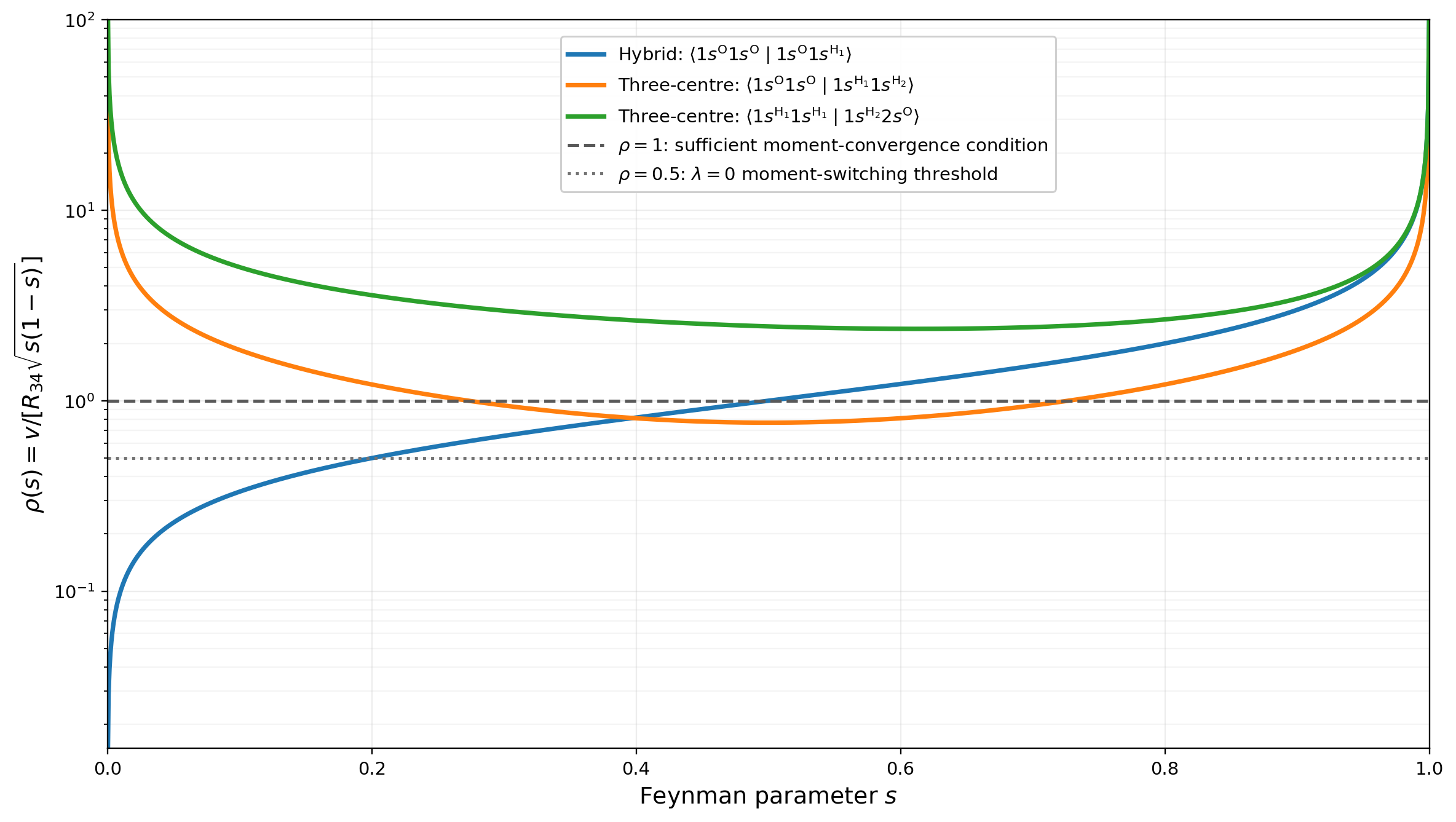}
\caption{Oscillation-to-damping ratio $\rho(s)$ for the semi-infinite integrals $\mathcal H(s)$ and $\mathcal K(s)$ arising in the H$_2$O hybrid and three-centre calculations. The horizontal axis is linear and the vertical axis is logarithmic (base 10). The lines $\rho=1$ and $\rho=0.5$ mark, respectively, the sufficient convergence condition for the moment expansion and the numerical switching threshold for $\lambda=0$.}
\label{rho-h2o-geometry}
\end{figure}

Figure~\ref{rho-h2o-geometry} shows the oscillation-to-damping ratio $\rho(s)$ for the semi-infinite integrals $\mathcal H(s)$ and $\mathcal K(s)$ associated with these H$_2$O configurations and a second three-centre assignment. The line $\rho=1$ marks the sufficient condition for convergence of the moment series, and $\rho=0.5$ indicates the switching threshold used in the $\lambda=0$ numerical implementation. The hybrid case enters this regime near $s=0$, whereas the three-centre case can enter it in the interior of the Feynman interval.

Table~\ref{small-v-diagnostics} gives the semi-infinite integrals evaluated by the moment expansion, their reference values, and the corresponding errors at selected nodes of the 48-point Gauss--Legendre quadrature for the outer Feynman integration. The reference values were obtained by direct integration of the complete original spherical-Bessel radial integrands using adaptive tanh--sinh quadrature, with a working precision of 45 decimal digits for both $\mathcal H(s)$ and $\mathcal K(s)$. The semi-infinite intervals were subdivided into finite subintervals and an infinite tail. These working precisions are not claimed certified error bounds.

Writing $\mathcal I(s)$ for either $\mathcal H(s)$ or $\mathcal K(s)$, we define the errors relative to the direct-integration reference values by:
\begin{equation}
\textrm{Error}_{\mathrm{mom}}=|\mathcal I_{\mathrm{mom}}(s)-\mathcal I_{\mathrm{ref}}(s)|,\qquad \textrm{Error}_{S+\mathrm{DE}}=|\mathcal I_{S+\mathrm{DE}}(s)-\mathcal I_{\mathrm{ref}}(s)|.
\label{small-v-errors}
\end{equation}

Table~\ref{small-v-diagnostics} reports the corresponding errors, with the moment expansion truncated at $q=Q$. For comparison, the $S+\mathrm{DE}$ quadrature is evaluated on a fixed grid with $M=16$, without switching to the moment expansion when $\rho(s)$ falls below the threshold used to activate this evaluation branch in the adaptive algorithm. This calculation serves as a diagnostic benchmark and is distinct from the adaptive implementation.

For the hybrid integral, the moment expansion agrees closely with direct integration near the lower endpoint, where $\rho(s)$ is small. For the three-centre H$_2$O configuration, $\rho(s)<1$ holds over an interior interval, with convergence slowing as $\rho(s)$ approaches unity. In the implemented algorithm, the moment expansion is automatically activated for $\rho\leq0.5$, with the switching threshold extended to $\rho\leq0.9$ when $\lambda\geq3$ in the $S+\mathrm{DE}$ method.

We next compare the two Fourier-type maps for the near-endpoint, strongly oscillatory kernel with $s=0.998$, $\lambda=4$, and $v=49.996$. Both maps use the condition $vMh=\pi$ and the same termwise truncation threshold of $10^{-15}$. Table~\ref{map-comparison} shows that both maps reach approximately the same accuracy, but $\phi_1$ requires substantially fewer quadrature nodes throughout the displayed refinement. We therefore use $\phi_1$ for the remaining $S+$DE calculations.

\begin{figure}[H]
\centering
\includegraphics[width=0.65\linewidth]{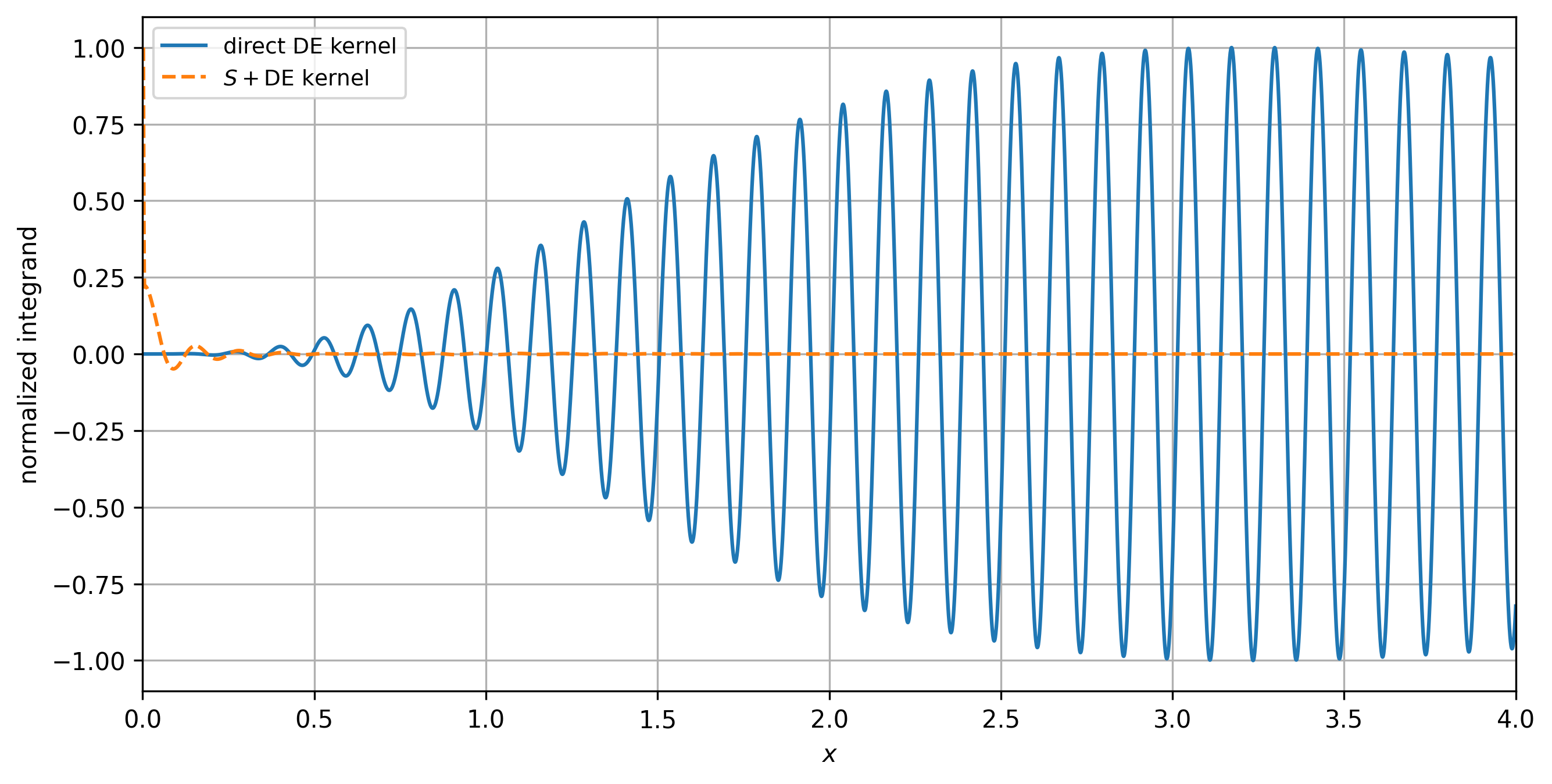}
\caption{Normalized oscillatory kernels before and after the $S$ transformation for the severe test case.}
\label{preintegrands}
\end{figure}

Figure~\ref{preintegrands} compares the oscillatory factors in a severe kernel. The $S$ transformation replaces the nonuniform finite zeros of $j_\lambda$ by the equally spaced sine zeros used in \eqref{auto-023}.

\begin{figure}[H]
\centering
{\small\textcolor{black}{$\bullet$}\ Direct DE\qquad \textcolor{blue}{$\bullet$}\ $S+\mathrm{DE}(\phi_1)$\qquad \textcolor{red}{$\bullet$}\ $S+\mathrm{DE}(\phi_2)$}\par\vspace{0.15ex}
{\scriptsize marker shape: $\bullet\ s=0.0024$, $\blacksquare\ s=0.5$, $\blacktriangle\ s=0.998$}\par\vspace{0.4ex}
\begin{tikzpicture}
\begin{loglogaxis}[width=.85\linewidth,height=3.35cm,xmin=8,xmax=650,ymin=60,ymax=20000,ylabel={nodes},title={${\lambda=2}$},title style={yshift=-0.6ex},grid=both,minor grid style={gray!15},major grid style={gray!30},xticklabels=\empty]
\addplot+[only marks,color=black,mark=*,mark size=3.0pt] coordinates {(102.1848,2150)(255.4619,5032)(510.9238,9587)};
\addplot+[only marks,color=black,mark=square*,mark size=3.0pt] coordinates {(10.0000,280)(25.0000,620)(50.0000,1165)};
\addplot+[only marks,color=black,mark=triangle*,mark size=3.0pt] coordinates {(111.9154,2344)(279.7884,5484)(559.5769,10453)};
\addplot+[only marks,color=blue,mark=*,mark size=3.0pt] coordinates {(102.1848,89)(255.4619,161)(510.9238,278)};
\addplot+[only marks,color=blue,mark=square*,mark size=3.0pt] coordinates {(10.0000,87)(25.0000,154)(50.0000,265)};
\addplot+[only marks,color=blue,mark=triangle*,mark size=3.0pt] coordinates {(111.9154,89)(279.7884,161)(559.5769,278)};
\addplot+[only marks,color=red,mark=*,mark size=3.0pt] coordinates {(102.1848,141)(255.4619,289)(510.9238,530)};
\addplot+[only marks,color=red,mark=square*,mark size=3.0pt] coordinates {(10.0000,133)(25.0000,270)(50.0000,491)};
\addplot+[only marks,color=red,mark=triangle*,mark size=3.0pt] coordinates {(111.9154,140)(279.7884,284)(559.5769,524)};
\end{loglogaxis}
\end{tikzpicture}
\vspace{-0.4ex}
\begin{tikzpicture}
\begin{loglogaxis}[width=.85\linewidth,height=3.35cm,xmin=8,xmax=650,ymin=60,ymax=20000,ylabel={nodes},title={${\lambda=4}$},title style={yshift=-0.6ex},grid=both,minor grid style={gray!15},major grid style={gray!30},xlabel={$\rho=v/[R_{34}\sqrt{s(1-s)}]$}]
\addplot+[only marks,color=black,mark=*,mark size=3.0pt] coordinates {(102.1848,2146)(255.4619,5024)(510.9238,9572)};
\addplot+[only marks,color=black,mark=square*,mark size=3.0pt] coordinates {(10.0000,276)(25.0000,612)(50.0000,1150)};
\addplot+[only marks,color=black,mark=triangle*,mark size=3.0pt] coordinates {(111.9154,2340)(279.7884,5476)(559.5769,10438)};
\addplot+[only marks,color=blue,mark=*,mark size=3.0pt] coordinates {(102.1848,81)(255.4619,108)(510.9238,176)};
\addplot+[only marks,color=blue,mark=square*,mark size=3.0pt] coordinates {(10.0000,81)(25.0000,106)(50.0000,173)};
\addplot+[only marks,color=blue,mark=triangle*,mark size=3.0pt] coordinates {(111.9154,81)(279.7884,108)(559.5769,176)};
\addplot+[only marks,color=red,mark=*,mark size=3.0pt] coordinates {(102.1848,93)(255.4619,175)(510.9238,305)};
\addplot+[only marks,color=red,mark=square*,mark size=3.0pt] coordinates {(10.0000,91)(25.0000,170)(50.0000,294)};
\addplot+[only marks,color=red,mark=triangle*,mark size=3.0pt] coordinates {(111.9154,93)(279.7884,175)(559.5769,305)};
\end{loglogaxis}
\end{tikzpicture}
\caption{Quadrature-node counts for direct DE, $S+\mathrm{DE}(\phi_1)$, and $S+\mathrm{DE}(\phi_2)$ over the complete nine-point $(s,v)$ grid for each $\lambda=2,4$.}
\label{regimemap}
\end{figure}
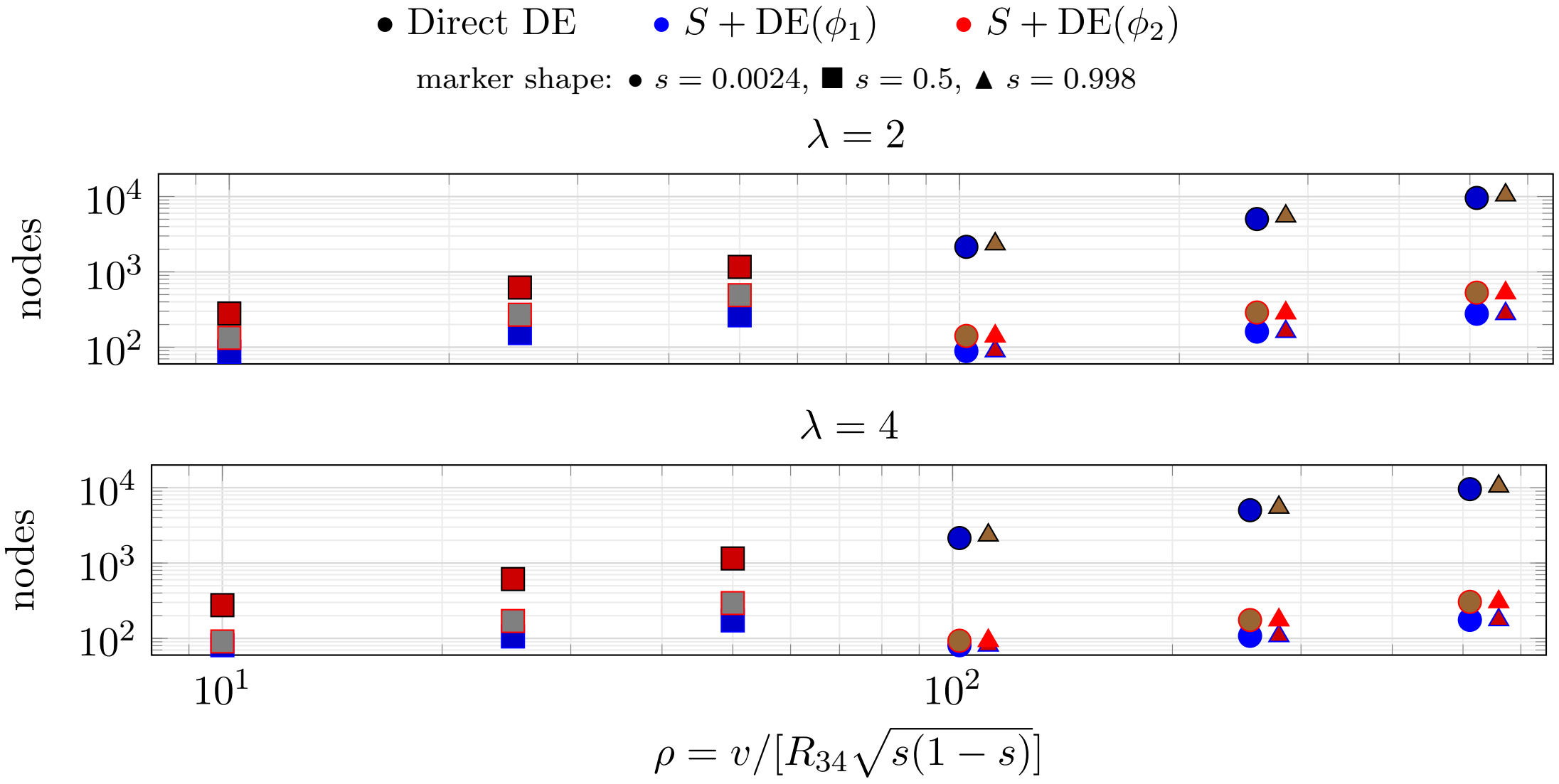

To compare the quadratures over different damping and oscillation regimes, we use the parameter $\rho(s)$ defined in \eqref{auto-017}. Figure~\ref{regimemap} presents the results as a function of $\rho(s)$~\eqref{auto-017} for $R_{34}=2$ and $\lambda=2$ and $4$, using the nine combinations of $s\in\{0.0024,\,0.5,\,0.998\}$ and $v\in\{10,\,25,\,50\}$. The quadrature parameters are $M=4$ and $h=\frac{\pi}{2vM}$ for direct DE, and $M=8,4,2$ with $h=\frac{\pi}{vM}$ for $S+\mathrm{DE}$. 

Figure~\ref{regimemap} compares direct DE with both $S+$DE mappings under these fixed quadrature parameters. The separation between direct DE and the two $S+$DE routes increases strongly with $\rho(s)$ and $\lambda$. For the near-endpoint, strongly oscillatory grid case $(s,v,\lambda)=(0.998,50,4)$, the direct DE and $S+\mathrm{DE}$ calculations require, respectively:
\begin{equation}
N_{\mathrm{direct\,DE}}=10438,\quad
N_{S+\mathrm{DE}(\phi_1)}=176
\quad\text{and}\quad
N_{S+\mathrm{DE}(\phi_2)}=305.
\label{severe-node-counts}
\end{equation}

To assess the sensitivity to the map parameter, we varied $K$ while holding $M$, the sine-zero alignment $h=\pi/(vM)$, and the stopping rule fixed in each case. Table~\ref{K-sensitivity} considers four radial kernels, specified by $(s,v,\lambda,M)$: Case I $(0.998,49.996,4,1)$, Case II $(0.5,10,2,8)$, Case III $(0.5,10,0,12)$, and Case IV $(0.0024,10,2,8)$. Thus the tests include both endpoint regions and the angular orders $\lambda=0,2,4$. Case III uses a finer fixed mesh to resolve the low-order interior kernel. For each case, the absolute errors are computed relative to the corresponding complete $S$-series value, held fixed as $K$ varies. At $K=6$, the errors range from approximately $3.7\times10^{-16}$ to $2.4\times10^{-15}$, with fewer nodes than for $K=3,4,5$. Larger values use fewer nodes but do not produce uniformly smaller errors. 

Table~\ref{severe-phase-resolution} gives the corresponding refinement comparison. In this table $h$ or $M$ is varied to show the convergence of direct DE and $S+$DE$(\phi_1)$ for the same kernel. The severe case reported in \cite{Safouhi19} clearly illustrates the effectiveness of the $S$ transformation. The direct-DE mesh with $h=0.04$ does not adequately resolve the oscillations, leading to an inaccurate evaluation of the integral. Reducing the spacing to $h=0.018$ requires 4695 nodes to reach an absolute error of $3.74\times10^{-15}$. By contrast, $S+$DE$(\phi_1)$ reaches an absolute error of $4.36\times10^{-16}$ with 107 nodes at $M=1$.

Table~\ref{convDE} shows the convergence with increasing $M$ for near-endpoint weak-damping and interior strong-damping cases. The interior case reaches essentially full retained accuracy at moderate $M$, while the near-endpoint cases require additional refinement as the molecular damping becomes weaker. For each calculation, the mesh spacing $h$ is determined from \eqref{auto-023}. 

Tables~\ref{h2o-threecentre} and \ref{h2o-hybrid} present the twenty H$_2$O three-centre and hybrid integrals reported in Table~9 of \cite{Safouhi19}. The results obtained with direct DE and $S+$DE$(\phi_1)$ are compared with the complete $S$-series and agree with the corresponding values in \cite{Safouhi19}. The complete molecular integrals were evaluated using a 48-point Gauss--Legendre quadrature for the outer Feynman-parameter integration.

The complete H$_2$O and C$_2$H$_4$ molecular integrals reported in Tables~\ref{h2o-threecentre}--\ref{hybrid-complete} are evaluated over normalized STFs. Each STF is first expanded into the finite linear combination of $B$ functions given in \eqref{stf-b-expansion}. Applying this expansion to all four orbital factors yields a finite sum of integrals over $B$ functions. Each such integral is evaluated using \eqref{auto-034}, with either direct DE or $S+$DE for its radial integral, and the results are combined with the coefficients of the STF--$B$ expansions. The tabulated values are compared with the complete $S$-series references and the available literature values \cite{Safouhi19}.

Tables~\ref{threecentre-complete} and \ref{hybrid-complete} give the extended C$_2$H$_4$ results with three-centre and hybrid integrals over STFs. In addition to the direct DE and $S+$DE$(\phi_1)$ evaluations, we report independently accumulated complete $S$-series values. The values reported in~\cite{Safouhi19} were compared with two independent STO codes, STOP~\cite{Bouferguene-Fares-Hoggan-57-801-96} and ADGGSTNGINT~\cite{Fernandez-Rico-Lopez-Ema-Ramirez-105-216-97}. Thus, agreement with these literature values also provides an indirect comparison with those earlier calculations.

\section{Numerical Tables and Figures}\label{numerical-data}

The tables and figures are presented in the order in which they are discussed in Section~\ref{numerical-analysis}.

\begin{table}[H]
\centering
\caption{Semi-infinite integrals $\mathcal H(s)$ and $\mathcal K(s)$ arising in the H$_2$O hybrid and three-centre calculations over $1s$ STFs.}
\label{small-v-diagnostics}
\setlength{\tabcolsep}{2.6pt}
\begin{tabular}{@{}lrrrrrrr@{}}
\toprule
Integral & $s$ & $\rho(s)$ & $Q$ & Reference & Moment & $\textrm{Error}_{\mathrm{mom}}$ & $\textrm{Error}_{S+\mathrm{DE}}$\\
\midrule
$\mathcal H(s)$ & 0.0006145 & 0.0248 & 6 & $1.083156042\times 10^{-4}$ & $1.083156042\times 10^{-4}$ & $9.49\times 10^{-19}$ & $1.36\times 10^{-2}$\\
$\mathcal H(s)$ & 0.0147042 & 0.1222 & 6 & $1.488961177\times 10^{-5}$ & $1.488961177\times 10^{-5}$ & $2.37\times 10^{-20}$ & $1.14\times 10^{-4}$\\
$\mathcal H(s)$ & 0.0617140 & 0.2565 & 6 & $7.073822747\times 10^{-7}$ & $7.073822747\times 10^{-7}$ & $4.53\times 10^{-18}$ & $1.00\times 10^{-6}$\\
$\mathcal H(s)$ & 0.1855663 & 0.4773 & 8 & $1.214983063\times 10^{-8}$ & $1.214983085\times 10^{-8}$ & $2.18\times 10^{-16}$ & $1.06\times 10^{-9}$\\
$\mathcal K(s)$ & 0.3256221 & 0.8993 & 20 & $5.510454669\times 10^{-6}$ & $5.510454970\times 10^{-6}$ & $3.02\times 10^{-13}$ & $4.74\times 10^{-11}$\\
$\mathcal K(s)$ & 0.4838099 & 0.7684 & 20 & $5.502668780\times 10^{-6}$ & $5.502668780\times 10^{-6}$ & $5.12\times 10^{-16}$ & $1.43\times 10^{-10}$\\
\bottomrule
\end{tabular}
\end{table}

\begin{table}[H]
\centering
\caption{Comparison of the two Fourier-type DE maps for the near-endpoint, strongly oscillatory kernel with $s=0.998$, $\lambda=4$, $v=49.996$, $R_{34}=2$, $\zeta_s=2$, and $\zeta_3=\zeta_4=1$. Both maps use $vMh=\pi$. The $S$-series reference value $\mathcal K_{\mathrm{ref}}=2.869930715021211\times10^{-4}$ and $\textrm{Error}=|\mathcal K-\mathcal K_{\mathrm{ref}}|$.} 
\label{map-comparison}
\begin{tabular}{lrrr}
\toprule
Method & $M$ & Nodes & \textrm{Error} \\
\midrule
$S+$DE$(\phi_1)$ & 1 & 107 & $4.36\times10^{-16}$ \\
$S+$DE$(\phi_1)$ & 2 & 176 & $5.94\times10^{-16}$ \\
$S+$DE$(\phi_1)$ & 4 & 308 & $9.53\times10^{-16}$ \\
\midrule
$S+$DE$(\phi_2)$ & 1 & 172 & $1.68\times10^{-15}$ \\
$S+$DE$(\phi_2)$ & 2 & 305 & $8.23\times10^{-16}$ \\
$S+$DE$(\phi_2)$ & 4 & 564 & $1.08\times10^{-15}$ \\
\bottomrule
\end{tabular}
\end{table}

\begin{table}[H]
\centering
\caption{Sensitivity of $S+$DE$(\phi_1)$ to the map parameter $K$ for four representative radial kernels (Cases I--IV, specified in the text). $N$ is the total quadrature-node count.}
\label{K-sensitivity}
\setlength{\tabcolsep}{3.1pt}
\begin{tabular}{@{}r *{4}{rr}@{}}
\toprule
& \multicolumn{2}{c}{Case I} & \multicolumn{2}{c}{Case II} & \multicolumn{2}{c}{Case III} & \multicolumn{2}{c}{Case IV} \\
\cmidrule(lr){2-3}\cmidrule(lr){4-5}\cmidrule(lr){6-7}\cmidrule(l){8-9}
$K$ & $N$ & Absolute Error & $N$ & Absolute Error& $N$ & $\textrm{Error}$ & $N$ & Absolute Error \\
\midrule
$3$ & 128 & $4.86\times10^{-16}$ & 158 & $5.38\times10^{-16}$ & 199 & $4.08\times10^{-14}$ & 174 & $2.13\times10^{-15}$ \\
$4$ & 117 & $1.06\times10^{-15}$ & 145 & $3.85\times10^{-16}$ & 182 & $3.02\times10^{-14}$ & 158 & $6.94\times10^{-16}$ \\
$5$ & 112 & $3.54\times10^{-16}$ & 136 & $8.00\times10^{-17}$ & 168 & $1.28\times10^{-14}$ & 148 & $1.10\times10^{-15}$ \\
$6$ & 107 & $4.36\times10^{-16}$ & 128 & $3.72\times10^{-16}$ & 158 & $2.37\times10^{-15}$ & 138 & $1.26\times10^{-15}$ \\
$2\pi$ & 105 & $7.46\times10^{-16}$ & 126 & $2.63\times10^{-16}$ & 155 & $2.59\times10^{-14}$ & 136 & $3.06\times10^{-15}$ \\
$8$ & 98 & $1.03\times10^{-15}$ & 116 & $1.19\times10^{-15}$ & 142 & $1.58\times10^{-14}$ & 123 & $1.09\times10^{-15}$ \\
\bottomrule
\end{tabular}
\end{table}

\begin{table}[H]
\centering
\caption{Comparison for the severe inner-kernel test of~\cite{Safouhi19}, with $s=0.998$, $\lambda=4$, $v=49.996$, $R_{34}=2$, $\zeta_s=2$, and $\zeta_3=\zeta_4=1$.}
\label{severe-phase-resolution}
\begin{tabular}{lrrr}
\toprule
Method & Control & Nodes & Absolute Error \\
\midrule
Direct DE & $M=4,h=0.040$ & 2189 & $1.21$ \\
Direct DE & $M=4,h=0.020$ & 4245 & $1.37\times10^{-13}$ \\
Direct DE & $M=4,h=0.018$ & 4695 & $3.74\times10^{-15}$ \\
\midrule
$S+$DE$(\phi_1)$ & $M=1$ & 107 & $4.36\times10^{-16}$ \\
$S+$DE$(\phi_1)$ & $M=2$ & 176 & $5.94\times10^{-16}$ \\
$S+$DE$(\phi_1)$ & $M=4$ & 308 & $9.53\times10^{-16}$ \\
\bottomrule
\end{tabular}
\end{table}

\begin{table}[H]
\centering
\caption{Refinement of the $S+$DE$(\phi_1)$ quadrature for representative inner radial kernels.}
\label{convDE}
\begin{tabular}{lrrr}
\toprule
$(s,v,\lambda)$ & $M$ & Nodes & Absolute Error\\
\midrule
$(0.0024,10,0)$ & 4 & 226 & $1.19\times10^{-11}$ \\
$(0.0024,10,0)$ & 8 & 237 & $1.07\times10^{-13}$ \\
\midrule
$(0.5,25,2)$ & 2 & 97 & $2.29\times10^{-15}$ \\
$(0.5,25,2)$ & 4 & 154 & $2.84\times10^{-17}$ \\
\midrule
$(0.998,49.996,4)$ & 1 & 107 & $4.36\times10^{-16}$ \\
$(0.998,49.996,4)$ & 2 & 176 & $5.94\times10^{-16}$ \\
\bottomrule
\end{tabular}
\end{table}

\begin{table}[H]
\centering
\caption{Complete three-centre two-electron Coulomb integrals over STFs for H$_2$O, with values reported in Table~9 of~\cite{Safouhi19}.}
\label{h2o-threecentre}
\setlength{\tabcolsep}{2.0pt}
\begin{tabular}{@{}p{4.0cm}rrr@{}}
\toprule
Integrals &  $S+$DE$(\phi_1)$ & Direct DE & $S$ series\\
\midrule
$\langle 1s^{\mathrm{O}}\,1s^{\mathrm{O}}\mid 1s^{\mathrm{H}_1}\,1s^{\mathrm{H}_2}\rangle$ & 0.166733422736855 & 0.166733422736836 & 0.166733422736899 \\
$\langle 1s^{\mathrm{O}}\,2s^{\mathrm{O}}\mid 1s^{\mathrm{H}_1}\,1s^{\mathrm{H}_2}\rangle$ & 0.033679727712632 & 0.033679727712653 & 0.033679727712789 \\
$\langle 1s^{\mathrm{O}}\,2p_{z}^{\mathrm{O}}\mid 1s^{\mathrm{H}_1}\,1s^{\mathrm{H}_2}\rangle$ & 0.001539018243091 & 0.001539018243118 & 0.001539018243109 \\
$\langle 2s^{\mathrm{O}}\,2s^{\mathrm{O}}\mid 1s^{\mathrm{H}_1}\,1s^{\mathrm{H}_2}\rangle$ & 0.145439173480290 & 0.145439173480310 & 0.145439173480327 \\
$\langle 2p_{z}^{\mathrm{O}}\,2p_{z}^{\mathrm{O}}\mid 1s^{\mathrm{H}_1}\,1s^{\mathrm{H}_2}\rangle$ & 0.134763477689259 & 0.134763477689266 & 0.134763477689288 \\
$\langle 2p_{-}^{\mathrm{O}}\,2p_{-}^{\mathrm{O}}\mid 1s^{\mathrm{H}_1}\,1s^{\mathrm{H}_2}\rangle$ & 0.128387519442823 & 0.128387519442819 & 0.128387519442835 \\
$\langle 2p_{+}^{\mathrm{O}}\,2p_{-}^{\mathrm{O}}\mid 1s^{\mathrm{H}_1}\,1s^{\mathrm{H}_2}\rangle$ & -0.001666457155435 & -0.001666457155432 & -0.001666457155437 \\
$\langle 1s^{\mathrm{H}_1}\,1s^{\mathrm{H}_1}\mid 1s^{\mathrm{H}_2}\,2s^{\mathrm{O}}\rangle$ & 0.218488995251726 & 0.218488995251714 & 0.218488995251758 \\
$\langle 1s^{\mathrm{H}_1}\,1s^{\mathrm{H}_1}\mid 1s^{\mathrm{H}_2}\,2p_{z}^{\mathrm{O}}\rangle$ & 0.157272683561797 & 0.157272683561803 & 0.157272683561822 \\
$\langle 1s^{\mathrm{H}_1}\,1s^{\mathrm{H}_1}\mid 1s^{\mathrm{H}_2}\,2p_{+}^{\mathrm{O}}\rangle$ & -0.084627891776168 & -0.084627891776161 & -0.084627891776178 \\
\bottomrule
\end{tabular}
\end{table}

\begin{table}[H]
\centering
\caption{Complete hybrid two-electron Coulomb integrals over STFs for H$_2$O, with values reported in Table~9 of~\cite{Safouhi19}.}
\label{h2o-hybrid}
\setlength{\tabcolsep}{2.0pt}
\begin{tabular}{@{}p{5.15cm}rrr@{}}
\toprule
Integrals & $S+$DE$(\phi_1)$ & Direct DE &  $S$ series\\
\midrule
$\langle 1s^{\mathrm{O}}\,1s^{\mathrm{O}}\mid 1s^{\mathrm{O}}\,1s^{\mathrm{H}_1}\rangle$ & 0.175630427352041 & 0.175630427352550 & 0.175630427352962 \\
$\langle 2s^{\mathrm{O}}\,1s^{\mathrm{O}}\mid 1s^{\mathrm{O}}\,1s^{\mathrm{H}_1}\rangle$ & 0.026546346644800 & 0.026546346649518 & 0.026546346649180 \\
$\langle 2p_{z}^{\mathrm{O}}\,1s^{\mathrm{O}}\mid 2p_{z}^{\mathrm{O}}\,1s^{\mathrm{H}_1}\rangle$ & 0.005284624781338 & 0.005284624781561 & 0.005284624781346 \\
$\langle 2p_{+}^{\mathrm{O}}\,1s^{\mathrm{O}}\mid 2p_{z}^{\mathrm{O}}\,1s^{\mathrm{H}_1}\rangle$ & 0.000883960239930 & 0.000883960240130 & 0.000883960239936 \\
$\langle 1s^{\mathrm{O}}\,1s^{\mathrm{O}}\mid 2p_{z}^{\mathrm{O}}\,1s^{\mathrm{H}_1}\rangle$ & 0.195970416002449 & 0.195970416002311 & 0.195970416002464 \\
$\langle 1s^{\mathrm{O}}\,1s^{\mathrm{O}}\mid 2p_{+}^{\mathrm{O}}\,1s^{\mathrm{H}_1}\rangle$ & 0.180590559485287 & 0.180590559485160 & 0.180590559485300 \\
$\langle 1s^{\mathrm{H}_1}\,1s^{\mathrm{H}_1}\mid 1s^{\mathrm{H}_1}\,1s^{\mathrm{O}}\rangle$ & 0.030786144383427 & 0.030786144383446 & 0.030786144383432 \\
$\langle 1s^{\mathrm{H}_1}\,1s^{\mathrm{H}_1}\mid 1s^{\mathrm{H}_1}\,2s^{\mathrm{O}}\rangle$ & 0.323596407152489 & 0.323596407152487 & 0.323596407173672 \\
$\langle 1s^{\mathrm{H}_1}\,1s^{\mathrm{H}_1}\mid 1s^{\mathrm{H}_1}\,2p_{z}^{\mathrm{O}}\rangle$ & 0.258864903094987 & 0.258864903095003 & 0.258864903135860 \\
$\langle 1s^{\mathrm{H}_1}\,1s^{\mathrm{H}_1}\mid 1s^{\mathrm{H}_1}\,2p_{+}^{\mathrm{O}}\rangle$ & -0.238549055692386 & -0.238549055692401 & -0.238549055730051 \\
\bottomrule
\end{tabular}
\end{table}

\begin{table}[H]
\centering
\caption{Complete three-centre two-electron Coulomb integrals over STFs for C$_2$H$_4$.}
\label{threecentre-complete}
\setlength{\tabcolsep}{1.6pt}
\begin{tabular}{@{}p{4.45cm}rrr@{}}
\toprule
Integrals & $S+$DE$(\phi_1)$ & Direct DE & Complete $S$ series\\
\midrule
$\langle 1s^{\mathrm{C}_1}\,1s^{\mathrm{C}_1}\mid 1s^{\mathrm{H}}\,1s^{\mathrm{C}_2}\rangle$ & 0.030219053645262 & 0.030219053646239 & 0.030219053645741 \\
$\langle 1s^{\mathrm{C}_1}\,1s^{\mathrm{C}_1}\mid 1s^{\mathrm{H}}\,2s^{\mathrm{C}_2}\rangle$ & 0.193446396651927 & 0.193446396652027 & 0.193446396651969 \\
$\langle 1s^{\mathrm{C}_1}\,1s^{\mathrm{C}_1}\mid 1s^{\mathrm{H}}\,2p_{z}^{\mathrm{C}_2}\rangle$ & 0.024517185798384 & 0.024517185799671 & 0.024517185798390 \\
$\langle 1s^{\mathrm{C}_1}\,1s^{\mathrm{C}_1}\mid 1s^{\mathrm{H}}\,2p_{+}^{\mathrm{C}_2}\rangle$ & 0.086095770850790 & 0.086095770851859 & 0.086095770850807 \\
$\langle 2s^{\mathrm{C}_1}\,1s^{\mathrm{C}_1}\mid 1s^{\mathrm{H}}\,1s^{\mathrm{C}_2}\rangle$ & 0.004958767237152 & 0.004958767237711 & 0.004958767237824 \\
$\langle 2s^{\mathrm{C}_1}\,2s^{\mathrm{C}_1}\mid 1s^{\mathrm{H}}\,1s^{\mathrm{C}_2}\rangle$ & 0.029165931997705 & 0.029165931997412 & 0.029165931997765 \\
$\langle 2p_{z}^{\mathrm{C}_1}\,2p_{z}^{\mathrm{C}_1}\mid 1s^{\mathrm{H}}\,2s^{\mathrm{C}_2}\rangle$ & 0.199079760183170 & 0.199079760183271 & 0.199079760183213 \\
$\langle 2p_{z}^{\mathrm{C}_1}\,2p_{z}^{\mathrm{C}_1}\mid 1s^{\mathrm{H}}\,2p_{z}^{\mathrm{C}_2}\rangle$ & 0.033248166777031 & 0.033248166776887 & 0.033248166777036 \\
\bottomrule
\end{tabular}
\end{table}

\begin{table}[H]
\centering
\caption{Complete hybrid two-electron Coulomb integrals over STFs for C$_2$H$_4$.}
\label{hybrid-complete}
\setlength{\tabcolsep}{1.6pt}
\begin{tabular}{@{}p{4.45cm}rrr@{}}
\toprule
Integrals & $S+$DE$(\phi_1)$ & Direct DE & Complete $S$ series\\
\midrule
$\langle 1s^{\mathrm{C}_1}\,1s^{\mathrm{C}_1}\mid 1s^{\mathrm{C}_1}\,2s^{\mathrm{C}_2}\rangle$ & 0.136049511034626 & 0.136049511034642 & 0.136049511034650 \\
$\langle 1s^{\mathrm{C}_1}\,1s^{\mathrm{C}_1}\mid 1s^{\mathrm{C}_1}\,2p_{z}^{\mathrm{C}_2}\rangle$ & -0.194543772192977 & -0.194543772194276 & -0.194543772194391 \\
$\langle 1s^{\mathrm{C}_1}\,1s^{\mathrm{C}_1}\mid 2p_{z}^{\mathrm{C}_1}\,2p_{z}^{\mathrm{C}_2}\rangle$ & -0.235395379625631 & -0.235395379625483 & -0.235395379625644 \\
$\langle 1s^{\mathrm{C}_1}\,1s^{\mathrm{C}_1}\mid 2p_{+}^{\mathrm{C}_1}\,2p_{+}^{\mathrm{C}_2}\rangle$ & 0.146236313435027 & 0.146236313435146 & 0.146236313435062 \\
$\langle 2p_{z}^{\mathrm{C}_1}\,2p_{z}^{\mathrm{C}_1}\mid 2p_{z}^{\mathrm{C}_1}\,2p_{z}^{\mathrm{C}_2}\rangle$ & -0.209119826213493 & -0.209119826208820 & -0.209119826208880 \\
$\langle 2p_{z}^{\mathrm{C}_1}\,2p_{z}^{\mathrm{C}_1}\mid 2p_{+}^{\mathrm{C}_1}\,2p_{+}^{\mathrm{C}_2}\rangle$ & 0.130266320162182 & 0.130266320159882 & 0.130266320159871 \\
$\langle 2p_{+}^{\mathrm{C}_1}\,2p_{+}^{\mathrm{C}_1}\mid 2p_{+}^{\mathrm{C}_1}\,2p_{+}^{\mathrm{C}_2}\rangle$ & 0.127717050483977 & 0.127717050485105 & 0.127717050485127 \\
\bottomrule
\end{tabular}
\end{table}

\section{Conclusions}\label{conclusions}

The radial kernels of three-centre two-electron Coulomb and hybrid integrals exhibit distinct geometry-dependent regimes. The ratio $\rho(s)$ of oscillation to molecular damping provides a sufficient convergence condition for the moment representation, $\rho(s)<1$. In hybrid integrals, this regime occurs naturally near the lower Feynman endpoint. In the case of three-centre integrals, it is determined by the centre geometry and can arise through interior cancellation of the displacement vectors.

Independent radial calculations show that the moment expansion can avoid large errors of an unmodified fixed-mesh $S+$DE rule in the small-$v$ regime. For larger oscillation-to-damping ratios, the $S$ transformation replaces the spherical-Bessel factor with a sine factor whose equally spaced zeros are adapted to Fourier-type double-exponential quadrature. The complete H$_2$O and C$_2$H$_4$ calculations over STFs provide molecular-integral values for comparison with the corresponding $S$-series evaluations and literature values. The radial comparisons justify the use of a small-$\rho$ evaluation branch, while certification of the switching thresholds at the complete molecular level remains a separate numerical question.

The combination of the $S$ transformation with DE quadrature provides a highly effective approach to the numerical evaluation of oscillatory molecular integrals, substantially reducing the required number of quadrature nodes while maintaining high accuracy. The complementary treatment of the weakly oscillatory regime through the convergent moment expansion further improves the reliability and versatility of the algorithm. Together, these developments represent an important step toward an efficient computational package for molecular electronic-structure calculations using exponential-type basis functions.

\section*{Statements and Declarations}
\noindent\textbf{Funding.} This research was supported by the Natural Sciences and Engineering Research Council of Canada (NSERC), Grant RGPIN-2016-04317 (extended to 2026).

\noindent\textbf{Competing Interests.} The authors declare that they have no competing interests.

\noindent\textbf{Data Availability.} No datasets were generated or analyzed during the current study.

\clearpage

\begin{thebibliography}{10}

\bibitem{Safouhi11}
H.~Safouhi.
\newblock An extremely efficient approach for accurate and rapid evaluation of three-centre two-electron {Coulomb} and hybrid integrals over {$B$} functions.
\newblock {\em J. Phys. A: Math. Gen.}, 34:881--902, 2001.

\bibitem{Safouhi19}
L.~Berlu and H.~Safouhi.
\newblock A new algorithm for accurate and fast numerical evaluation of hybrid and three-center two-electron {Coulomb} integrals over {Slater} type functions.
\newblock {\em J. Phys. A: Math. Gen.}, 36:11267--11283, 2003.

\bibitem{Safouhi60}
R.M. Slevinsky and H.~Safouhi.
\newblock Compact formulae for three-center nuclear attraction integrals over exponential type functions.
\newblock {\em Journal of Mathematical Chemistry}, 60:1337--1355, 2022.

\bibitem{Safouhi38}
H.~Safouhi.
\newblock {Bessel}, sine and cosine functions and extrapolation methods for computing molecular multi-center integrals.
\newblock {\em Numer. Algor.}, 54:141--167, 2010.

\bibitem{Safouhi18}
L.~Berlu and H.~Safouhi.
\newblock An extremely efficient and rapid algorithm for a numerical evaluation of three-center nuclear attraction integrals over {Slater} type functions.
\newblock {\em J. Phys. A: Math. Gen.}, 36:11791--11805, 2003.

\bibitem{Safouhi20}
L.~Berlu and H.~Safouhi.
\newblock Multicenter two-electron {Coulomb} and exchange integrals over {Slater} functions evaluated using a generalized algorithm based on nonlinear transformations.
\newblock {\em J. Phys. A: Math. Gen.}, 37:3393--3410, 2004.

\bibitem{Safouhi21}
H.~Safouhi.
\newblock Highly accurate numerical results for three-center nuclear attraction and two-electron {Coulomb} and exchange integrals over {Slater} type functions.
\newblock {\em Int. J. Quantum Chem.}, 100(2):172--183, 2004.

\bibitem{Lesiuk-Moszynski-90-063318-14}
M.~Lesiuk and R.~Moszynski.
\newblock Reexamination of the calculation of two-center, two-electron integrals over {Slater}-type orbitals. i. {Coulomb} and hybrid integrals.
\newblock {\em Physical Review E}, 90:063318, 2014.

\bibitem{Yu-Si-Zhou-123-e27128-23}
T.~Yu, Y.~Si, and J.~Zhou.
\newblock General expressions of molecular multicenter integrals for slater-type orbitals.
\newblock {\em International Journal of Quantum Chemistry}, 123:e27128, 2023.

\bibitem{Stark-Meier-Hatch-Kammeraad-Dang-Zimmerman-47-e70291-26}
A.~Stark, N.~Meier, J.~Hatch, J.~A. Kammeraad, D.-K. Dang, and P.~M. Zimmerman.
\newblock Numerical integration of {Slater} basis functions over prolate spheroidal grids.
\newblock {\em Journal of Computational Chemistry}, 47:e70291, 2026.

\bibitem{Kaur-Chen-Luhana-Silberman-Ewert-Bunger-Segal-Greif-47-e70472-26}
A.~Kaur, W.~Chen, G.~K. Luhana, L.~Silberman, M.~Ewert, A.~B{\"u}nger, U.~Segal, and C.~Greif.
\newblock Fast and accurate gpu-accelerated computation of two-electron four-center coulomb repulsion integrals.
\newblock {\em Journal of Computational Chemistry}, 47:e70472, 2026.

\bibitem{Safouhi55}
J.~Lovrod and H.~Safouhi.
\newblock Double exponential transformation for computing three-center nuclear attraction integrals.
\newblock {\em Molecular Physics}, 118(4):e1619854, 2020.

\bibitem{Steinborn-Filter-38-273-75}
E.O. Steinborn and E.~Filter.
\newblock Translations of fields represented by spherical-harmonics expansions for molecular calculations. {III}. {Translations} of reduced {Bessel} functions, {Slater}-type s-orbitals, and other functions.
\newblock {\em Theor. Chim. Acta.}, 38:273--281, 1975.

\bibitem{Filter-Steinborn-18-1-78}
E.~Filter and E.O. Steinborn.
\newblock Extremely compact formulas for molecular two-center one-electron integrals and {Coulomb} integrals over {Slater}-type atomic orbitals.
\newblock {\em Phys. Rev. A.}, 18(1):1--11, 1978.

\bibitem{Condon-Shortley-51}
E.U. Condon and G.H. Shortley.
\newblock {\em The theory of atomic spectra}.
\newblock Cambridge University Press, Cambridge, England, 1951.

\bibitem{Watson-44}
G.N. Watson.
\newblock {\em A Treatise on the Theory of {Bessel} Functions}.
\newblock Cambridge University Press, Second Edition, Cambridge, England, 1944.

\bibitem{Abramowitz-Stegun-65}
M.~Abramowitz and I.A. Stegun.
\newblock {\em Handbook of Mathematical Functions}.
\newblock Dover, New York, 1965.

\bibitem{DLMF-26}
{NIST Digital Library of Mathematical Functions}.
\newblock {NIST Digital Library of Mathematical Functions}, 2026.
\newblock Chapters 3, 10, and 13; sections used across the seven manuscripts include 3.6, 10.22, 10.29, 10.32, 10.40, 10.43, 10.51, 10.74, and 13.6; release 1.2.7; accessed August 2026.

\bibitem{Slater-42-33-32}
J.C. Slater.
\newblock Analytic atomic wave functions.
\newblock {\em Phys. Rev.}, 42(1):33--43, 1932.

\bibitem{Safouhi15}
H.~Safouhi.
\newblock Convergence properties of the {$S\overline{D}$} transformation and a fast and accurate numerical evaluation of molecular integrals.
\newblock {\em J. Phys. A: Math. Gen.}, 35:9685--9698, 2002.

\bibitem{Denich-Novati-198-160-24}
E.~Denich and P.~Novati.
\newblock Some notes on the trapezoidal rule for {Fourier} type integrals.
\newblock {\em Applied Numerical Mathematics}, 198:160--175, 2024.

\bibitem{Denich-Novati-100-23-24}
E.~Denich and P.~Novati.
\newblock A {Sinc} rule for the {Hankel} transform.
\newblock {\em Journal of Scientific Computing}, 100:23, 2024.

\bibitem{Dominguez-26}
V.~Dom\'inguez.
\newblock Fast and stable computation of highly oscillatory and/or exponentially decaying integrals using a {Clenshaw--Curtis} product-integration rule.
\newblock {\em Numerical Algorithms}, pages 1--40, 2026.
\newblock Online publication; volume and final pagination not yet assigned in the bibliographic record consulted.

\bibitem{Takahasi-Mori-9-721-74}
H.~Takahasi and M.~Mori.
\newblock Double exponential formulas for numerical integration.
\newblock {\em RIMS}, 9(3):721--741, 1974.

\bibitem{Mori-12-119-85}
M.~Mori.
\newblock Quadrature formulas obtained by variable transformation and the {DE}-rule.
\newblock {\em Journal of Computational and Applied Mathematics}, 12--13:119--130, 1985.

\bibitem{Ooura-Mori-38-353-91}
T.~Ooura and M.~Mori.
\newblock The double exponential formula for oscillatory functions over the half infinite interval.
\newblock {\em J. Comput. Appl. Math.}, 38:353--360, 1991.

\bibitem{Ooura-Mori-112-229-99}
T.~Ooura and M.~Mori.
\newblock A robust double exponential formula for fourier-type integrals.
\newblock {\em J. Comput. Appl. Math.}, 112:229--241, 1999.

\bibitem{Mori-Sugihara-127-287-01}
M.~Mori and M.~Sugihara.
\newblock The double-exponential transformation in numerical analysis.
\newblock {\em Journal of Computational and Applied Mathematics}, 127(1,2):287--296, 2001.

\bibitem{Mori-41-897-05}
M.~Mori.
\newblock Discovery of the double exponential transformation and its developments.
\newblock {\em RIMS}, 41(4):897--935, 2005.

\bibitem{Sugihara-72-767-02}
M.~Sugihara.
\newblock Near optimality of the {Sinc} approximation.
\newblock {\em Mathematics of Computation}, 72(242):767--786, 2003.

\bibitem{Bouferguene-Fares-Hoggan-57-801-96}
A.~Bouferguene, M.~Fares, and P.~E. Hoggan.
\newblock {STOP}: A {Slater}-type orbital package for molecular electronic structure determination.
\newblock {\em International Journal of Quantum Chemistry}, 57(4):801--810, 1996.

\bibitem{Fernandez-Rico-Lopez-Ema-Ramirez-105-216-97}
J.~Fern\'{a}ndez Rico, R.~L\'{o}pez, I.~Ema, and G.~Ram\'{i}rez.
\newblock Calculation of many-centre two-electron molecular integrals with {STO}.
\newblock {\em Computer Physics Communications}, 105:216--224, 1997.

\end{thebibliography}

\end{document}